\documentclass[12pt]{amsart}
\usepackage{amsmath,amsfonts,amssymb,amscd, epsfig}
\usepackage{latexsym}
\usepackage{amsthm}

\usepackage[%dvipdfm, 
hyperindex=true,pagebackref=true,bookmarks=true,
colorlinks=true,linkcolor=blue,citecolor=red]
{hyperref}
\def\~{{\rm --}} 
\title [Affine non-gatherable triples] 
{Non-gatherable triples for classical affine root systems}
\author[Ivan Cherednik]{Ivan Cherednik $^\dag$}
\author[Keith Schneider]{Keith Schneider}
\thanks{$^\dag$  \today\ \ \ Partially supported by NSF grant
DMS--0800642}

\address[I. Cherednik]{Department of Mathematics, UNC
Chapel Hill, North Carolina 27599, USA\\
chered@email.unc.edu}
\address[K. Schneider]{Department of Mathematics, UNC
Chapel Hill, North Carolina 27599, USA\\
schneidk@email.unc.edu}

 \def\bysame{{\bf --- }}
 \def\~{{\bf --}}
\newcommand{\comment}[1]{}
\renewcommand{\tilde}{\widetilde}
\renewcommand{\hat}{\widehat}
\renewcommand{\bar}{\overline}
\def\si{\hbox{\em \ae}}

\renewcommand{\tilde}{\widetilde}
\renewcommand{\hat}{\widehat}

\newcommand{\Z}{{\mathbb Z}}

\newcommand{\R}{{\mathbb R}}

\def\HH{\mbox{${\mathcal H}$\kern-5.2pt${\mathcal H}$}}

\newtheorem{theorem}{Theorem}[section]

\newtheorem{proposition}[theorem]{Proposition}

\newtheorem{lemma}[theorem]{Lemma}
\newtheorem{corollary}[theorem]{Corollary}

\newtheorem{theorem }{Theorem}[section]
\newtheorem{maintheorem }[theorem]{Main Theorem}
\newtheorem{proposition }[theorem]{Proposition}
\newtheorem{definition }[theorem]{Definition}
\newtheorem{lemma }[theorem]{Lemma}
\newtheorem{corollary }[theorem]{Corollary}
\newtheorem{notation }[theorem]{Notation}
\newtheorem{remark }[theorem]{Remark}
\newtheorem{example }[theorem]{Example}

\newtheorem{ maintheorem }[theorem]{Main Theorem}
\newtheorem{ theorem}{Theorem}[section]
\newtheorem{ proposition}[theorem]{Proposition}
\newtheorem{ definition}[theorem]{Definition}
\newtheorem{ lemma}[theorem]{Lemma}
\newtheorem{ corollary}[theorem]{Corollary}
\newtheorem{ notation}[theorem]{Notation}
\newtheorem{ remark}[theorem]{Remark}
\newtheorem{ example}[theorem]{Example}

 \newcommand{\rmk}{{\bf Comment.\ }}

\def\for{\  \hbox{ for } \ }
\def\iif{ \ \hbox{ if } \ }

\def\where{\  \hbox{ where } \ }
\def\and{\  \hbox{ and } \ }
\def\and{\  \hbox{ and } \ }

\def\equal{\stackrel{\,\mathbf{def}}{= \kern-3pt =}}

\def\la{\lambda}
\def\La{\Lambda}
\def\om{\omega}

\def\th{\theta}
\def\al{\alpha}
\def\be{\beta}
\def\ga{\gamma}
\def\ep{\epsilon}

\def\de{\delta}

\def\si{\sigma}

\def\Ga{\Gamma}
\def\ze{\zeta}

\def\vth{{\vartheta}}

\def\tal{\tilde{\alpha}}
\def\tbe{\tilde{\beta}}

\def\tla{\tilde{\lambda}}
\def\tga{\tilde{\gamma}}
\def\tGa{\tilde{\Gamma}}

\def\tu{\tilde{u}}

\def\tw{\widetilde w}
\def\tW{\widetilde W}

\def\tz{\tilde z}

\def\tR{\tilde R}

\def\hw{\widehat{w}}
\def\hW{\widehat{W}}
\def\hu{\hat{u}}

\def\0{\mathbf{0}}

\def\çF{\mathcal{F}}

\def\lan{\langle}

\def\ran{\rangle}

\newcommand{\sq}{\phantom{1}\hfill$\qed$}

\def\HH{\mathfrak{H}}

\def\HH{\hbox{${\mathcal H}$\kern-5.2pt${\mathcal H}$}}

\font\smm=msbm10 at 12pt 
\def\symbol#1{\hbox{\smm #1}}
\def\lsmash{{\symbol n}}

\def\#{\sharp}

\begin{document}

\comment{
This paper contains a complete description of minimal non-gatherable 
triangle triples in the lambda-sequences for the affine
classical root systems and some claims for arbitrary
(reduced) affine root systems. It continues our previous paper
devoted to the non-affine case; interestingly, the affine theory
clarifies the classification in the non-affine case.  
The lambda-sequences are associated with reduced decompositions 
(words) in affine Weyl groups. The existence of the 
non-gatherable triples is a combinatorial obstacle for using the 
technique of intertwiners in the theory of 
irreducible representations of the (double) affine Hecke algebras, 
complementary to their algebraic-geometric theory.}

\maketitle
%\comment{  %%%%% PART I
\renewcommand{\baselinestretch}{1.0} 
{\small
\tableofcontents
}          %%%%% PART I
\renewcommand{\baselinestretch}{1.0}
%} 
%$\mathfrak{a,b,c,d,e,f,g,h,i,j,k,l,m,n,o,p,q,r,s,t,e^u,e^v,e^w}$
\vfill

\renewcommand{\natural}{\wr}

\setcounter{section}{-1}
\setcounter{equation}{0}
\section{Introduction}

This paper is a continuation of \cite{CS} and the
part of \cite{C0} devoted to 
{\em non-gatherable triangle triples} in $\la$\~sequences. 
The $\la$\~sequences are the 
sequences of positive roots associated with 
reduced decompositions (words) in affine and nonaffine
Weyl groups. The minimal {\em non-gatherable triangle triples},
NGT, $\{\al,\al+\be,\be\}$
is a $\la$\~sequences with non-movable 
(under the Coxeter 
transformations) endpoints $\al,\be$ such that
$\al+\be$ is a root and $|\al|=|\be|$.
Their nonaffine classification
for $B_n,C_n (n\ge 3), D_n (n\ge 4)$ and for $F_4,E_6$
is the subject of \cite{CS}; 
there are no NGT for nonaffine and affine
$A_n, B_2, C_2, G_2$.

We describe all 
minimal NGT for the classical affine root systems based
on their planar interpretation from \cite{Ch5} and provide 
a universal general construction for arbitrary (reduced,
irreducible) root systems. In principle, the latter can
be used to obtain all such triples. 
{\em A posteriori}, only a significantly reduced version
of our universal theorem is sufficient for the 
{\em classical} affine NGT, however, it can be more
involved for the exceptional root systems.
\smallskip

The affine minimal NGT we construct are given in terms 
{\em almost dominant} weights (where one simple root
can be disregarded in the definition of the dominant
weights) and certain ``small"
elements from the nonaffine Weyl group. 
The weight itself generally is not sufficient to determine the 
corresponding minimal affine NGT uniquely. 

For the classical
affine root systems, the answer appeared very explicit. 
Combinatorially, it is given in terms of partitions
of a type $A$ subdiagram inside the initial 
nonaffine Dynkin diagram and (additionally)
an increasing sequences of non-negative
integers associated with such partitions. 
\smallskip

Interestingly, all such minimal NGT can be naturally 
presented in terms of those of type $B$. It is not 
unexpected because the planar interpretation unifies all 
classical root systems in one construction.
The passage to the other types from $B$ is by using certain 
{\em parity corrections} directly related to the element  
$s_0$ of type $B$ treated as an element of the {\em extended}
affine Weyl groups of type $C,D$ (in the twisted setting).  
\smallskip

We note that there is a natural general (all root systems)
procedure for producing candidates for {\em nonaffine} 
minimal NGT  from the affine ones. For the exceptional 
root systems,
it is justified only modulo certain technical assumptions. 
Employing it, we come to a somewhat surprising construction
of {\em nonaffine} minimal NGT in affine terms. It 
clarifies our nonaffine classification for the classical root 
systems and gives a promising approach to the exceptional 
nonaffine root systems, where the description of all minimal 
NGT is known but is technically involved.   
\smallskip 

The existence of NGT is a combinatorial obstacle for 
using the technique of intertwiners (see, e.g. \cite{C0})
in the theory of irreducible representations
of the affine and double affine Hecke algebras, complementary 
to the geometric approach from of \cite{KL1} and its double
affine generalization. We mainly mean the constructive
theory of such representations (where the intertwining 
elements are used to construct basic vectors).

The theory of affine and double affine algebras motivated
our paper a great deal, but NGT are quite interesting in their 
own right.
Gathering together the triangle triples
using the Coxeter transformations seems an important question 
in the theory of reduced decompositions of Weyl groups, 
which is far from being simple. 
More generally, assuming that
$\la(w)$ contains all positive roots of a certain root subsystem, 
can these roots be gathered using the Coxeter transformations?
\medskip

{\bf Basic definitions.}
Let $R\in \R^n$ be a reduced irreducible root system or
its affine extension, $W$ the corresponding Weyl group.
Then the $\la$\~{\em set} is defined as 
$\la(w)=R_+\cap w^{-1}(-R_+)$ for $w\in W$, where
$R_+$ is the set of positive roots in $R$. It is well-known
that $w$ is uniquely determined by  
$\la(w)$; many properties of $w$ and its
reduced decompositions can be interpreted in terms
of this set. The $\la$\~{\em sequence\,}
is the $\la$\~set with the ordering of roots naturally
induced by a given reduced decomposition.  

The intrinsic description of such sets and sequences
is mainly given in terms of the {\em triangle triples} 
$\{\be,\ga=\al+\be,\al\}$. For instance, $\al,\be\in \la(w)$
$\Rightarrow$ $\al+\be\in \la(w)$ and the latter 
root must appear between $\al$ and $\be$\, if this set
is treated as a sequence. This property is necessary but
not sufficient; see \cite{C0} for a 
comprehensive discussion.
 
We want to know when the sets of positive 
roots of rank two subsystems inside a
given {\em sequence} 
$\la(w)$ can be {\em gathered} (made consecutive) 
using the Coxeter transformations in $\la(w)$. It is natural
to allow the transformations only within
the minimal segments containing these roots.
This problem can be readily 
reduced to considering the {\em triangle triples}
provided some special conditions on the lengths.
The answer is always affirmative
only for the root systems $A_n, B_2, C_2, G_2$ (and their
affine counterparts) or in the case when $|\al|\neq|\be|$. 
Otherwise non-trivial NGT always exist.
\medskip

{\bf The planar representation.}
For the root system $A_n$ (nonaffine or affine), 
gathering the triples is simple. It readily results from the 
planar interpretation of 
the reduced decompositions and the corresponding
$\la$\~sequences in terms of $(n+1)$ lines in the 
two-dimensional plane (on the cylinder in the affine
case). 

Conceptually, this interpretation is 
a variant of the classical {\em geometric} approach to the 
reduced decompositions of $w\in W$ in terms of the lines 
(or pseudo-lines) that go 
from the main Weyl chamber to the chamber
corresponding to $w$; see \cite{Bo}. However, 
the planar description adds a lot to this general 
approach. It is a powerful tool, which
dramatically simplifies dealing with combinatorial
problems concerning the reduced decompositions.   
\smallskip

The $A_n$\~planar interpretation  was extended in \cite{Ch0}
to other {\em classical} root systems and $G_2$,
and then to their affine extensions in \cite{Ch5}. 
Omitting $G_2$,
it is given in terms of $n$ lines in $\R^2$ with reflections 
in one {\em mirror} for the nonaffine
$B_n,C_n,D_n$ and two {\em mirrors} in the affine case.
This approach is significantly developed in this paper;
it can be used for quite a few problems beyond NGT. 
 
We were able to use the planar interpretation  
to find {\em all} minimal non-gatherable triples, 
{\em minimal NGT}, for the affine root systems $B,C,D$. 
Algebraically,
without such geometric support, it is an involved 
combinatorial problem. 
No planar (or similar) interpretation is known for 
$F_4, E_{6,7,8}$. Nonaffine minimal NGT can be 
classified  using computers (see \cite{CS} for $F_4,E_6$);
the exceptional affine root systems will be considered in  
our further works. 
\smallskip

Generally, the {\em admissibility} condition from \cite{C0} is 
necessary and sufficient for the triple to be {\em gatherable},
which is formulated in terms of subsystems of $R$ of types 
$B_3,C_3$ or $D_4$. This universal (but not very convenient
to use) theorem can be now re-established for the classical
root systems using the classification we give in this paper. 
\medskip

{\bf Relation to (double) affine Hecke algebras.}
The existence of NGT and some other features of similar 
nature are not present in the case of $A$. Generally, the theory
of root systems is uniform at level of generators and 
relations of the corresponding Weyl (or braid) groups;
however the root systems behave quite differently when the 
``relations for Coxeter relations" are considered.
\smallskip

Presumably, the phenomenon of NGT is one of the major 
combinatorial obstacles for creating a universal
theory of AHA-DAHA ``highest vectors"
generalizing Zelevinsky's {\em segments} in the $A$\~case
and based on the intertwining operators.
This technique was fully developed
only for affine and double affine Hecke algebras
of type $A_n$ and in some cases of small ranks.

The classification and explicit description
of {\em semisimple\,} irreducible representations of AHA and DAHA 
is expected to be a natural application of this technique.
The recent research (in progress) indicates that a thorough 
analysis of NGT is needed for this and similar projects.
\smallskip

The fact that all triples are 
gatherable in the case of $A_n$ was the key in
\cite{Ch1}  and quite a few
further papers on the {\em quantum fusion procedure\,}.
This procedure reflects the duality of
AHA and DAHA of type $A$ are the corresponding
quantum groups and quantum toroidal algebras.  

Quantum groups and Yangians certainly deserve
special comments. In the case of $GL$,
their irreducible representations 
can be described in terms of the
so-called fusion procedure. The key object of the
latter is the {\em transfer matrix}, a
product of quantum $R$\~matrices 
geometrically corresponding to a bunch 
of $n$ parallel lines 
intersecting another bunch of $m$ parallel lines.

Major parts of this big theory were extended to 
the $R$\~matrices {\em with reflection} and
the {\em twisted Yangians} (of reflection type). The
the corresponding {\em transfer matrices}
are associated with the following configurations.
The $n$\~bunch of lines intersects 
the $m$\~bunch parallel to the mirror, then
reflects in this mirror and then again 
intersects the $m$\~bunch.
There are interesting modifications here when $D$ is
considered. These configurations (when $n\ge 2$)
are exactly those for the {\em non-affine} minimal 
NGT of type $B,C$. 
Recent research on the twisted Yangians \cite{KN}
indicates that  it is not by chance and that 
minimal NGT may be of importance for this theory.  

Expanding the theory of transfer matrices  to the 
affine case is a natural challenge, including
the {\em corner transfer matrices}, which are also related 
to our constructions. We hope that the 
classification of classical affine minimal NGT configurations 
will play its role in the (future) theory of twisted 
Yangians and Quantum groups of toroidal type.
\medskip

\section{Affine Weyl groups}
\setcounter{equation}{0}
Let $R=\{\al\}   \subset \R^n$ be a root system of type
$A,B,...,F,G$
with respect to a Euclidean form $(z,z')$ on $\R^n
\ni z,z'$,
$W$ the {\em Weyl group} \index{Weyl group $W$}
generated by the reflections $s_\al$,
$R_{+}$ the set of positive  roots ($R_-=-R_+$)
corresponding to fixed simple 
roots $\al_1,...,\al_n,$
$\Ga$ the Dynkin diagram
with $\{\al_i, 1 \le i \le n\}$ as the vertices.

We will also use sometimes
the dual roots (coroots) and the dual root system:
$$R^\vee=\{\al^\vee =2\al/(\al,\al)\}.$$

The root lattice and the weight lattice are:
\begin{align}
& Q=\oplus^n_{i=1}\Z \al_i \subset P=\oplus^n_{i=1}\Z \om_i,
\notag
\end{align}
where $\{\om_i\}$ are fundamental weights:
$ (\om_i,\al_j^\vee)=\de_{ij}$ for the
simple coroots $\al_i^\vee.$
Replacing $\Z$ by $\Z_{\pm}=\{m\in\Z, \pm m\ge 0\}$ we obtain
$Q_\pm, P_\pm.$ Here and further see  \cite{Bo}.
%Note that $Q\cap P_+\subset Q_+.$ It readily gives that
%each $\om_j$ has all
%strictly positive coefficients (generally, rational) 
%when expressed in terms of
%$\{\al_i\}.$

The form will be normalized
by the condition  $(\al,\al)=2$ for 
{\em short} roots. 
When dealing with the classical root
systems, the {\em most natural inner product $(\,,\,)_\ep$}
is the one making the $\ep_i$ in \cite{Bo} orthonormal.
It coincides with our $(\,,\,)$ for $C$ and $D$;
in the case of $B$, our form is $2(\,,\,)_\ep$.
One has:

\centerline{
$\nu_\al\equal (\al,\al)/2$ can be either $1,$ or $\{1,2\},$ or
$\{1,3\}.$ }

This normalization leads to the inclusions
$Q\subset Q^\vee,  P\subset P^\vee,$ where $P^\vee$ is
defined to be generated by
the fundamental coweights $\om_i^\vee.$

\smallskip
Let  $\vth\in R^\vee $ be the {\em maximal positive
coroot}. Equivalently, it is {\em maximal positive
short root} in $R$ due to our choice of the normalization.
All simple roots appear in its decomposition in $R$
or $R^\vee$.
Note that $2\ge (\vth,\al^\vee)\ge 0$ for $\al>0,$
$(\vth,\al^\vee)=2$ only for $\al=\vth,$ and
$s_{\vth}(\al)<0$ if $(\vth,\al)>0.$

%\vskip 0.2cm
\subsection{Affine roots}
The vectors $\ \tal=[\al,\nu_\al j] \in
\R^n\times \R \subset \R^{n+1}$
for $\al \in R, j \in \Z $ form the
\index{affine root system} {\em affine root system}
$\tR \supset R$ ($z\in \R^n$ are identified with $ [z,0]$).
We add $\al_0 \equal [-\vth,1]$ to the simple
 roots for the
{\em maximal short root}
\index{maximal short root $\vth$} $\vth$.
The corresponding set
$\tR_+$ of positive roots coincides with %%%% tR_+BOOK!!!
$R_+\cup \{[\al,\nu_\al j],\ \al\in R, \ j > 0\}$.

We will write $\tR=\widetilde{A}_n,
\widetilde{B}_n,\ldots, \widetilde{G}_2$ when dealing
with classical root systems.

The root system $\tR_+$ is called the 
{\em twisted} affine extension of $R$. The standard
one from \cite{Bo} is defined for maximal {\em long} root
$\th\in R_+$ and with omitting $\nu_\al$ in the
expression for the affine roots; the
inner product is normalized by the condition $(\th,\th)=2$.
The transformation of our considerations to the non-twisted 
case is straightforward.

Any positive affine root $[\al,\nu_\al j]$
is a linear combinations
with non-negative integral coefficients of
$\{\al_i,\,0\le i\le n\}$.
Indeed, it is well known that 
$[\al^\vee,j]$ is such combination
in terms of  $\{\al_i^\vee,\, 1\le i\le n\}$
and $[-\vth,1]$ for the system of affine {\em coroots},
that is $\tR^\vee=\{ [\al^\vee,j],\, \al\in R,\, j\in\Z \}$.
Hence, $[-\al,\nu_\al j]$  $=\nu_\al[-\al^\vee,j]$ has
the required representation.

Note that the sum of the long roots is always long,
the sum of two short roots can be a long root only
if they are orthogonal to each other. 

We complete the Dynkin diagram $\Ga$ of $R$
by $\al_0$ (by $-\vth$, to be more
exact); it is called {\em affine Dynkin diagram}
$\tGa$. One can obtain it from the
completed (extended by zero) Dynkin diagram 
from \cite{Bo} for the {\em dual system}
$R^\vee$ by reversing all arrows. 
%The number of laces between $\al_i$ and
%$\al_j$ in $\tGa$ will be denoted by $m_{ij}.$

The set of
the indices of the images of $\al_0$ by all
the automorphisms of $\tGa$ will be denoted by $O$
($O=\{0\} \for E_8,F_4,G_2$). Let $O'=\{r\in O, r\neq 0\}$.
The elements $\om_r$ for $r\in O'$ are the so-called minuscule
weights: $(\om_r,\al^\vee)\le 1$ for
$\al \in R_+$.

Given $\tal=[\al,\nu_\al j]\in \tR,  \ b \in P$, let
\begin{align}
&s_{\tal}(\tz)\ =\  \tz-(z,\al^\vee)\tal,\
\ b'(\tz)\ =\ [z,\ze-(z,b)]
\label{ondon}
\end{align}
for $\tz=[z,\ze] \in \R^{n+1}$.

The 
{\em affine Weyl group} $\tW$ is generated by all $s_{\tal}$
(we write $\tW = \lan s_{\tal}, \tal\in \tR_+\ran)$. One can take
the simple reflections $s_i=s_{\al_i}\ (0 \le i \le n)$
as its
generators and introduce the corresponding notion of the
length. This group is
the semidirect product $W\lsmash Q'$ of
its subgroups $W=$ $\lan s_\al,
\al \in R_+\ran$ and $Q'=\{a', a\in Q\}$, where
\begin{align}
& \al'=\ s_{\al}s_{[\al,\,\nu_{\al}]}=\
s_{[-\al,\,\nu_\al]}s_{\al}\for
\al\in R.
\label{ondtwo}
\end{align}

The
{\em extended Weyl group} $ \hW$ generated by $W\and P'$
(instead of $Q'$) is isomorphic to $W\lsmash P'$:
\begin{align}
&(wb')([z,\ze])\ =\ [w(z),\ze-(z,b)] \for w\in W, b\in P.
\label{ondthr}
\end{align}
{\em From now on,  $b$ and $b',$ $P$ and $P'$ 
will be identified.}

Note that the extended affine Weyl group in the standard
(non-twisted case) is identified with the semidirect product 
$W\lsmash P^\vee$.

The action in $\R^{n+1}$ is dual to the {\em affine action} 
$\hw(\!(z)\!)\equal w(z+\xi b)$ in $\R^n\ni z$ for a free 
parameter $\xi$, where  $\hw=wb$ and $w\in W, b\in P$ . 
I.e., $P$ acts via the 
translations in this definition. In more detail,
let $(\,[z,t]\,,\,z'\,)_\xi\equal
(z,z')+\xi t$ 
For  $z,z'\in \R^n,\ t\in \R$ and $\hw=wb\in \hW$,
\begin{align}\label{standact}
&(\,\hw(z)\,,\,\hw(\!(z')\!)\,)_\xi\ =\ 
(\,z\,,\,z'\,)_\xi.
\end{align}
Note that $s_{[\al,j]}(\!(z)\!)=z-2((z,\al)+j\xi)\al^\vee$

Given $b\in P_+$, let $w^b_0$ be the longest element
in the subgroup $W_0^{b}\subset W$ of the elements
preserving $b$. This subgroup is generated by simple
reflections. We set
\begin{align}
&u_{b} = w_0w^b_0  \in  W,\ \pi_{b} =
b( u_{b})^{-1}
\ \in \ \hW, \  u_i= u_{\om_i},\pi_i=\pi_{\om_i},
\label{xwo}
\end{align}
where $w_0$ is the longest element in $W,$
$1\le i\le n.$

The elements $\pi_r\equal\pi_{\om_r}, r \in O'$ and
$\pi_0=\hbox{id}$ leave $\tGa$ invariant
and form a group denoted by $\Pi$,
 which is isomorphic to $P/Q$ by the natural
projection $\{\om_r \mapsto \pi_r\}$. As to $\{ u_r\}$,
they preserve the set $\{-\vth,\al_i, i>0\}$.
The relations $\pi_r(\al_0)= \al_r= ( u_r)^{-1}(-\vth)$
distinguish the
indices $r \in O'$. Moreover (see e.g., \cite{C0}):
\begin{align}
& \hW  = \Pi \lsmash \tW, \where
  \pi_rs_i\pi_r^{-1}  =  s_j \iif \pi_r(\al_i)=\al_j,\
 0\le j\le n.
\end{align}
\medskip
%\vskip 0.2cm
\subsection{The length}
Setting
$\hw = \pi_r\tw \in \hW$ for $\pi_r\in \Pi, \tw\in \tW,$
the length $l(\hw)$
is by definition the length of the reduced decomposition
$\tw= $ $s_{i_l}...s_{i_2} s_{i_1} $
in terms of the simple reflections
$s_i, 0\le i\le n.$ 

The {\em length} can be
also defined as the
cardinality $|\la(\hw)|$
of the {\em $\la$\~set} of $\hw$\,:
\begin{align}\label{lasetdef}
&\la(\hw)\equal\tR_+\cap \hw^{-1}(\tR_-)=\{\tal\in \tR_+,\
\hw(\tal)\in \tR_-\},\
\hw\in \hW.
\end{align}

Note that $\la(\hw)$ is closed with respect to positive linear
combinations. More exactly, if $\tal=u\tbe+v\tga\in \tR$ for
rational
$u,v>0$, then $\tal\in \la(\hw)$ if  $\tbe\in \la(\hw) \ni \tga$.
Vice versa, if $\la(\hw)\ni \tal=u\tbe+v\tga$ for
$\tbe,\tga\in \tR_+$ and rational $u,v>0$, then either $\tbe$ or
$\tga$
must belong to $\la(\hw)$. Also,
\begin{align}\label{laclosed}
&\tal=[\al,\nu_\al j]\in\la(\hw)\, \Rightarrow\,
[\al,\nu_\al i]\in \la(\hw)\notag\\
&\hbox{for\ all\ } 0\le i < j\hbox{\ where\ }
i>0 \hbox{\ as\ } \al<0\,.
\end{align}
\smallskip

The coincidence with the previous definition
is directly related to the equivalence of
the following four claims:
\begin{align}
&(a)\ \ \  l(\hw\hu)=
l(\hw)+l(\hu)
\for \hw,\hu\in\hW\ (\hbox{length formula}),
\label{ltutwa}\\
&(b)\ \ \  \la(\hw\hu) = \la(\hu) \cup
\hu^{-1}(\la(\hw))\ (\hbox{cocycle relation}),
\label{ltutw}\\
&(c)\ \ \   \hu^{-1}(\la(\hw))
\subset \tR_+ \ (\hbox{positivity condition}),
\label{ltutwc}\\
&(d)\ \ \  \la(\hu)\subset \la_\nu(\hw)\ 
(\hbox{embedding condition}).
\label{ltutwd}
\end{align}
The key here is the following general relation:
\begin{align}
&  \la(\hw\hu)\ =\ \la(\hu)\,\widetilde{\cup}\,\,
\hu^{-1}(\la(\hw)) \hbox{ \ for \ any \ }
\hu,\hw,
\label{latutw}
\end{align}
where, by definition, the {\em reduced union}
$\,\widetilde{\cup}\,$ is obtained from\, $\cap$\, 
upon the cancelation of all pairs $\{\tal,-\tal\}$.
In particular, (\ref{latutw}) gives that
$$
\la(\hw^{-1})\ =\ -\hw(\la(\hw)).
$$

Applying (\ref{ltutw}) to the reduced decomposition
$\hw=\pi_rs_{i_l}\cdots s_{i_2}s_{i_1},$
\begin{align}
\la(\hw) = \{\ &\tal^l=\tw^{-1}s_{i_l}(\al_{i_l}),\
\ldots,\ \tal^3=s_{i_1}s_{i_2}(\al_{i_3}),\notag\\
&\tal^2=s_{i_1}(\al_{i_2}),\ \tal^1=\al_{i_1}\  \}.
\label{tal}
\end{align}
It demonstrates directly that the cardinality
$l$ of the set $\la(\hw)$ equals $l(\hw).$
Cf. \cite{Hu},4.5.

\rmk
It is worth mentioning that counterparts of the
$\la$\~sets can be introduced for 
$w=s_{i_l}\cdots s_{i_2}s_{i_1}$
in arbitrary Coxeter groups. Following \cite{Bo}
(Ch. IV, 1.4, Lemma 2), one can define 
\begin{align}
\La(w) = \{\ &t_l=w^{-1}s_{i_l}(s_{i_l}),\
\ldots,\ t_3=s_{i_1}s_{i_2}(s_{i_3}),\notag\\
&t_2=s_{i_1}(s_{i_2}),\ t_1=s_{i_1}\  \},
\label{talc}
\end{align}
where the action is by conjugation; $\La(w)\subset W$.

The $t$\~elements are (all) pairwise different if and
only if the decomposition is reduced (a simple 
straight calculation; see \cite{Bo}). 
Then this set does not depend
on the choice of the reduced decomposition. It readily 
gives a proof of formula (\ref{tal}) by induction 
and establishes the equivalence
of (a),(b) and (c). 

Generally,
the crystallographical case is significantly simpler
than the case of abstract Coxeter groups; using the root 
systems  dramatically simplifies
theoretical and practical (via computers) analysis of 
the reduced decompositions.  
The positivity of roots, 
the alternative definition of the $\la$\~sets from 
(\ref{lasetdef}) and, more specifically, 
property (c) are (generally)
missing in the theory of abstract Coxeter groups.
\sq
\medskip

In this paper, we will mainly treat $\la(\hw)$ as
sequences, called $\la$\~{\em sequences};
the roots in (\ref{tal}) are ordered naturally.
The sequence structures of the same $\la$\~set 
correspond to different choices of the reduced
decompositions of $\hw$.

An arbitrary simple root $\al_i\in\la(\hw)$ can be made
the first in a certain $\la$\~sequence. More generally:
\begin{align}
&\la(\hw)\
=\ \{\al>0, \ l( \hw s_{\al}) \le l(\hw) \};
\label{xlambda1}
\end{align}
see \cite{Bo} 
and \cite{Hu},4.6, Exchange Condition.

\smallskip
The {\em sequence} $\la(\tw)=\{\tal^l,\ldots,\tal^1\}$, where
$l=l(\tw)$, determines $\tw\in \tW$ uniquely. Indeed,
\begin{align} \label{lambdainv}
&\al_{i_1}=\tal^1,\, \al_{i_2}=s^1(\tal^2), \ldots,
\al_{i_p}=s^1s^2\cdots s^{p-1}(\tal^p),\ldots\notag\\
&\al_{i_l}=s^1s^2\cdots s^{l-1}(\tal^l)\,,\where \notag\\
&s^p=s_{\tal^p},\ \hw=s_{i_l}\cdots s_{i_1}=
s^1\cdots s^l.
\end{align}
Notice the order of the reflections $s^p$ 
in the decomposition of $\tw$ is {\em inverse}.
Moreover,
$\la(\hw)$ considered as an {\em unordered set}
determines $\hw$ uniquely up to the left multiplication by   
the elements $\pi_r\in \Pi$.
\smallskip

The intrinsic definition of the $\la$\~sequences is as follows.

$(i)$ Assuming that $\tal,\tbe,\tga=\tal+\tbe\in \tR_+$,
if $\tal,\tbe\in \la$ then $\tga\in\la$ and $\tga$
appears between $\tal,\tbe$;
if $\tal\not\in\la$ then $\tbe$ belongs to $\la$
and appears in $\la$ before $\tga$.

$(ii)$ If $\tal=[\al,\nu_\al j]\in \la$ then
$[\al,\nu_\al j\,'\,]\in \la$ as $j>j\,'>0$ and it appears in
$\la$ before $\tal$.

$(iii)$ If $\tbe\in \la$ and $\tga=\tbe-[\al,\nu_\al j]\in 
\tR_+[-]$
for $\al\in R_+,\,j\ge 0$, then $\tga\in \la$ and it
appears before $\tbe$.
\smallskip

If $\la$ is treated as an unordered set, then
it is in the form 
$\la=\la(\hw)$ for some $\hw\in \hW$ if and only if 
($i+ii+iii$) are imposed without the claims concerning
the ordering.

\medskip

\medskip
\subsection{Reduction modulo 
\texorpdfstring{\mathversion{bold}$W$}{\em W}}
It generalizes the construction of the elements
$\pi_{b}$ for $b\in P_+.$ 

\begin{proposition} \label{PIOM}
 Given $ b\in P$, there exists a unique decomposition
$b= \pi_b  u_b,$
$ u_b \in W$ satisfying one of the following equivalent conditions:

{(i) \  } $l(\pi_b)+l( u_b)\ =\ l(b)$ and
$l( u_b)$ is the greatest possible,

{(ii)\  }
$ \la(\pi_b)\cap R\ =\ \emptyset$.

The latter condition implies that
$l(\pi_b)+l(w)\ =\ l(\pi_b w)$
for any $w\in W.$ Besides, the relation $ u_b(b)
\equal b_-\in P_-=-P_+$
holds, which, in its turn,
determines $ u_b$ uniquely if one of the following equivalent
conditions is imposed:

{(iii) }
$l( u_b)$ is the smallest possible,

{(iv)\ }
if\, $\al\in \la( u_b)$  then $(\al,b)\neq 0$.
\end{proposition}
\qed

Condition (ii) readily
gives a complete description of the $\la$\~sets corresponding
to the elements $\pi_P=\{\pi_b, b\in P\}$. The roots there
must be from 
$$
\la\subset \tR_+[-]\equal\,\{[\al,\,\nu_\al j], \al\in R_-,j>0\}.
$$
only
See (\ref{laclosed}) and
(\ref{lambpi}) below. 
\smallskip

Setting $\tal=[\al,\nu_\al j]\in \tR_+,$ one has:
\begin{align}
\la(b) = \{ \tal,\  &( b, \al^\vee )>j\ge 0 \iif \al\in R_+,
\label{xlambi}\\
&( b, \al^\vee )\ge j> 0 \iif \al\in R_-\},
\notag \\
\la(\pi_b) = \{ \tal,\ \al\in R_-,\
&( b_-, \al^\vee )>j> 0
\iif  u_b^{-1}(\al)\in R_+,
\label{lambpi} \\
&( b_-, \al^\vee )\ge j > 0 \iif
u_b^{-1}(\al)\in R_- \}, \notag \\
\la(\pi_b^{-1}) = \{ \tal\in \tR_+,\  -&(b,\al^\vee)>j\ge 0 \},
\label{lapimin}\\
\la(u_b)\ =\ \{ \al\in R_+,\ \ \ \,&(b,\al^\vee)> 0 \}.
\label{laumin}
\end{align}

\smallskip
The element $b_{-}= u_b(b)$ is a unique element
from $P_{-}$ that belongs to the orbit $W(b)$.
Thus the equality   $c_-=b_- $ means that $b,c$
belong to the same orbit. We will also use
$b_{+} \equal w_0(b_-),$ a unique element in $W(b)\cap P_{+}.$
In terms of the elements $\pi_b,$
$$u_b\pi_b\ =\ b_-,\ \pi_b^{-1} u_b^{-1}\ =\ \varsigma(b_+),\ 
\varsigma(b)\equal -w_0(b).$$

Note that $l(\pi_b w)=l(\pi_b)+l(w)$ for all $b\in P,\ w\in W.$
For instance,
\begin{align}
&l(b_- w)=l(b_-)+l(w),\ l(wb_+)=l(b_+)+l(w).
\label{lupiw}
\end{align}
\smallskip

The definition of $\pi_b$ and $u_b$ is compatible
with the one from (\ref{xwo}) when
$b\in P_+$. Namely, 
\begin{align}
&u_{b} = w_0w^b_0  \in  W,\ \pi_{b} =
b( u_{b})^{-1}
\ \in \ \hW \for b\in P_+.
\label{xwox}
\end{align}
Recall that $w^b_0$ is the longest element
in the subgroup $W_0^{b}\subset W$ of the elements
preserving $b$, $w_0$ is the longest element in $W.$

We will need below this construction extended to
arbitrary $b\in P$ as follows. Let $w^b_0$ be the longest element
in the subgroup $W_0^{b}\subset W$, defined as the span
of {\em simple} reflections $s_i(1\le i\le n)$ preserving $b$. 
We set
\begin{align}
&\upsilon_{b}\equal w_0w^b_0  \in  W,\ \varpi_{b} =
b( \upsilon_{b})^{-1}
\ \in \ \hW.
\label{xwoy}
\end{align}
For $b\in P_+$,
the group $W_0^{b}$ coincides with the complete centralizer
of $b$ in $W$;\ $\upsilon_b=u_b$ and $\varpi_b=\pi_b$. 
Note that $w_0w^b_0=w^{b^\varsigma}_0w_0$ and 
\begin{align}
&\varpi_b^{-1}=\varpi_{\varsigma(b)}
\for \varsigma(b)=-w_0(b).
\label{xwoinverse}
\end{align}

\setcounter{equation}{0}
\section{General theory of NGT}
The transformations of the reduced decompositions
in $\tW$ are generated by the elementary ones, the
{\em Coxeter transformations}, that are substitutions
$(\cdots s_is_js_i)\mapsto (\cdots s_js_is_j)$ in
reduced decompositions of the elements $\tw\in \tW$. The number
of $s$\~factors is $2,3,4,6$ when $\al_i$ and $\al_j$ are connected
by $0,1,2,3$ laces in the affine or nonaffine 
Dynkin diagram. These transformations induce
{\em reversing the order} of the corresponding segments
(with $2,3,4,6$ roots)
of $\la(\tw)$ treated as a sequence. 
These segments can be naturally identified 
with the standard sequences of positive roots 
of type $A_1\times A_1$, $A_2$, $B_2$ or $G_2$.
The conjugations by $\pi_r\in \Pi$ will be applied too;
they permute of the indices of the
words from $\tW$ (preserving the length). 
\smallskip

\subsection{Admissibility condition}
The theorem below is essentially from \cite{C0}; it 
has application to the decomposition of the
polynomial representation of DAHA and is important
for the classification
of semisimple representations of AHA and DAHA (in
progress). We think that it clarifies
why dealing with the intertwining
operators for arbitrary root systems is significantly more 
difficult than in the $A_n$\~case (where much is known).
\smallskip

%From now on we will consider only $\tW$;
%the passage to $\hW$ does not add anything new as far as 
%NGT are concerned.
Given a reduced decomposition of $\hw\in \hW$,
let us assume that $\tal+\tbe=\tga$ for the roots
$\ldots,\tbe,\ldots,\tga,\ldots,\tal\ldots$
in $\la(\hw)$ ($\tal$ appears the first), where
only the following combinations of their lengths
are allowed in the $\widetilde{B},
\widetilde{C},\widetilde{F}$ cases
\begin{align}
&\hbox{long}+\hbox{long}=\hbox{long}\ \,
 (B,F_4) \hbox{\, \ or\ \,}
\hbox{short}+\hbox{short}=\hbox{short} \ \,
(C,F_4).
\label{shtshtsht}
\end{align}
We call such $\{\tbe,\tga,\tal\}$ a (triangle)
{\em triple}.

Since we will use the Coxeter transformations only
inside the
segment $[\tbe,\tal]\subset\la(w)$, from $\tal$ to
$\tbe$, it suffices to assume that  $\tal$ is a simple root.
The root systems 
$\widetilde{A}_n,
\widetilde{B}_2, \widetilde{C}_2, \widetilde{G}_2$ are excluded
from the following theorem; there are no NGT in these cases.

\begin{theorem}\label{RANKTWO}
The roots $\tbe,\tga,\tal$ from a triple are non-gatherable, i.e.,  
cannot be made consecutive roots
using the Coxeter transformations inside the segment
$[\tbe,\tal]\subset\la(\hw)$ if and only if a root
subsystem of type 
$B_3$, $C_3$ or $D_4$  exists such that
its intersection with $\la(\hw)$ 
constitutes the $\la$\~set of a certain
non-gatherable triple there. \sq
\end{theorem}
\medskip

The theorem can be readily reduced to considering
the elements $\hw$ representing 
{\em minimal NGT}, i.e., such that the $\la$\~sequence
$\la(\hw)$ begins with $\tal$ and ends with $\tbe$
and both roots (the endpoints) are non-movable
with respect to the Coxeter transformations of $\hw$.
Thus a minimal NGT is a pair,
the triple $\{\tbe,\tga=\tal+\tbe,\tal\}$ and the element
$\hw\in \hW$ that {\em represents} this triple.
Since such triple is uniquely determined by
$\hw$, we will constantly call $\hw$  minimal NGT
too, somewhat abusing the terminology.

The classification of the classical affine minimal NGT 
gives this statement for 
the {\em classical} affine root systems.
For the exceptional 
root systems, the first (universal, for all root
systems) part of the paper can be used. 
\smallskip

\subsection{Almost dominant weights}
Recall that we defined 
\begin{align}
&\upsilon_{b}\equal w_0w^b_0  \in  W,\ \varpi_{b} =
b\,(\upsilon_{b})^{-1}
\ \in \ \hW
\label{xwoyy}
\end{align}
for an arbitrary $b\in P$ (not only for dominant
ones). Given $b$, 
let us remove from the (nonaffine)
Dynkin diagram $\Ga$ the
vertices $\al_j$ such that $((\al_j,b)\neq 0$ and
represent the output 
as a union of connected subdiagrams $\Ga^{(m)}$. Thus:
$$
\Ga\ =\ \cup_m\Ga^{(m)}\cup\{\al_j\}\ \hbox {such that}\ 
\{\al_j\}=\{\al_j\mid (b,\al_j)\neq 0\}.
$$   
We will denote $\cup_m\Ga^{(m)}$ by $\Ga^b$;
$w^b_0$ is the product $\prod_m w_0^{(m)}$ for
the longest elements in the Weyl groups $W^{(m)}$
defined for $\Ga^{(m)}$.
Note that $\Ga^{\varsigma(b)}=
\varsigma(\Ga^b)$ for $\varsigma(b)=b^\varsigma=-w_0(b).$

This definition will be mainly used for 
``almost dominant" $b$ (when one simple root
is omitted). Let us fix a nonaffine simple
root $\al_k$ ($1\le k\le n$). The weight $b\in P$ in the 
constructions below will be always assumed from
\begin{align} \label{Pplusk}
&P_+^{(k)}\equal\{\,a\in P\, \mid\, (a,\al_j)\ge 0 \for 
j\ne k,\, (a,\al_k)\le 0\,\}.
\end{align}

Notice that we allow here $(b,\al_k)=0$. Then 
the corresponding $b$ will be dominant. In this case,
there can be several choices for $k$; we pick one $k$
such that $(b,\al_k)=0$, 
and construct $\dot{w}_0^b$ for 
$\Ga^b_\cdot\equal\Ga^b\setminus \{\al_k\}$.  Thus,  
$\dot{w}_0^b$
depends on the choice of $k$ in this case. Accordingly, 
let $R_{+\cdot}^{b}$ be the set of positive roots of 
the root system associated with $\Ga^{b}_\cdot$ 
considered as a subsystem
of $\Ga$. We will use below that
$\la(\dot{w}_0^b)=R_{+\cdot}^{b}$.

Let 
$\upsilon_{b}^\cdot=w_0\dot{w}_0^b$ and
$\varpi_{b}^\cdot=b(\upsilon_b^\cdot)^{-1}$.
For the sake of
uniformity, the {\em dot}\~notation 
will be used when $(b,\al_k)<0$; no 
{\em dot}\~modifications of $\Ga^b$, $w_0^b, \upsilon_b$
and $\varpi_b$ are 
necessary in this case. 

We need to extend the construction of $\varpi_{b}^\cdot$ even 
further by allowing certain reductions, the {\em $\si$\~reductions}, 
of the elements $\upsilon_b^\cdot$. For $\si\in W$, we set
\begin{align}
&\upsilon^{\si}_{b}\ =\ \varsigma(\si)\upsilon_b^\cdot= 
w_0\si\dot{w}_0^b 
\ \in\  W,\ \, \varpi^\si_{b}\ =\ 
b\,(\upsilon^\sigma_{b})^{-1}
\ \in \ \hW,\notag\\ 
&(\varpi_b^{\si})^{-1}=
\varpi_{\bar{b}}^{\bar{\si}} \for
\bar{b}\equal  \varsigma(\si(b)), \ 
\bar{\si}\equal \dot{w}_0^{\varsigma({b})}\, 
\varsigma(\si^{-1})\,
\dot{w}_0^{\bar{b}}.\label{barbsi} 
\end{align}
Accordingly,
\begin{align*}
&\upsilon^{\bar{\si}}_{\bar{b}} = (\upsilon^{\si}_{b})^{-1}
=w_0\bar{\si}\dot{w}_0^{\bar{b}} 
\in  W,\ \, \varpi^{\bar{\si}}_{\bar{b}} = 
\bar{b}\,(\upsilon^{\bar{\si}}_{\bar{b}})^{-1}.
\end{align*}
The expression for $\bar{\si}$ is equivalent to
the following relation:
\begin{align}\label{barsidotw}
&\bar{\si}\dot{w}_0^{\bar{b}}\ =\ 
\dot{w}_0^{\varsigma(b)} \varsigma(\si^{-1})\ =\ 
\varsigma(\si\dot{w}_0^{\varsigma(b)})^{-1}.
\end{align}
We note that $-\upsilon^{\si}_b(b)=\varsigma(\si(b))=\bar{b}$,
which readily gives that $\varsigma(\bar{\si}(\bar{b}))=b$.
Thus, our {\em bar}\~operation is involutive (by construction).

Recall that the definitions of
$\dot{w}_0^b, \dot{w}_0^{\bar{b}}$ 
depend on the choice of $k,\bar{k}$ when $(\al_k,b)=0$
and $(\al_{\bar{k}},\bar{b})=0$ (they can be not unique such).
The notation $\upsilon_b^{\si}$ automatically includes
the {\em dot}\~extension;
if $\si=\,$id then $\upsilon_b^{\si}=\upsilon_b^\cdot$
and $\varpi_b^{\si}=\varpi_b^\cdot$.

Generally, the decompositions in (\ref{barsidotw})
are not reduced. Let us address it.

\begin{proposition}\label{BARBLEN}
For $b\in P_+^{(k)}$, 
we take $\si\in W$ such that
$\bar{b}\equal \varsigma(\si(b))\in P_+^{(\bar{k})}$
for certain $1\le \bar{k}\le n$. Then the
following conditions
\begin{align}\label{lengdotwnew}
&(a):\ \ l(\bar{\si} \dot{w}_0^{\bar{b}})=
l(\bar{\si})+l(\dot{w}_0^{\bar{b}}),\ \  
&(b):\ \ \ l(\si \dot{w}_0^{b})=
l(\si)+l(\dot{w}_0^{b}).
\end{align}
are correspondingly equivalent to:
\begin{align*}
&(\tilde{a}):\, \al_j\not\in\la(\varpi^\si_b)\ \hbox{for}\
0\neq j\neq \bar{k},\ \ (\tilde{b}):\, 
\al_j\not\in\la(\varpi^{\bar{\si}}_{\bar{b}})\ \hbox{for}\
0\neq j\neq k.
\end{align*}
%If $(\al_{\bar{k}},\bar{b})=0$ or $(\al_k,b)=0$ then
%$\al_{\bar{k}}\not\in\la(\varpi^\si_b)$ or 
%$\al_{k}\not\in\la(\varpi^{\bar{\si}}_{\bar{b}})$
%correspondingly.
\end{proposition} 
{\em Proof}. It suffices to consider $(a)$;
the case of $(b)$ is analogous (and formally
follows from $(a)$). One has: 
\begin{align}\label{varpibsi}
\varpi_b^\si\ =\ b (w_0\si\dot{w}_0^{b})^{-1}=
(\dot{w}_0^{b}\si^{-1} w_0)\cdot(-\bar{b})=
(w_0\bar{\si}\dot{w}_0^{\bar{b}})\cdot(-\bar{b}).
\end{align}
The set $\la(-\bar{b})$ contains $\al\in R_+$
if and only if $(-\bar{b},\al)>0$
due to formula (\ref{xlambi}). Therefore,  
this set contains a simple
nonaffine root only when $(-\bar{b},\al_{\bar{k}})>0$;
then it can be only $\al_{\bar{k}}$.
We use that $\bar{b}\in P_{+\cdot}^{(\bar{k})}$.
Other nonaffine simple roots 
$\la(\varpi_b^\si)$ can come only from
$\bar{b}(\la((w_0\bar{\si}\dot{w}_0^{\bar{b}}))$.

Conditions $(a,b)$ are equivalent to:
\begin{align}\label{upsilonsharp}
&(a'):\ \, \la(\bar{\si} \dot{w}_0^{\bar{b}})\ =\ 
R_{+\cdot}^{\,\bar{b}}\,\cup\, 
\dot{w}_0^{\bar{b}}(\la(\bar{\si})),\\ 
&(b'):\ \, \la(\si \dot{w}_0^{b})\ =\ 
R_{+\cdot}^{\,b}\,\cup\, 
\dot{w}_0^{b}(\la(\si)).
\end{align}
Recall that  $\la(\dot{w}_0^{\bar{b}})=R_{+\cdot}^{\bar{b}}$\,,
where the latter is the set of all positive
roots in the subsystem with simple (nonaffine) roots
$\al_j$ such that $(\al_j,\bar{b})=0$ subject to the
following {\em dot}\~modification. 
If $(\al_{\bar{k}},\bar{b})=0$
(which is allowed), then $\al_{\bar{k}}$ must be
excluded from $\{\al_j\}$. 

We will use that
$\la(w_0 u)=R_+\setminus \la(u)$ for
any $u\in W$,
which is obvious from the definition of the $\la$\~sets. 

Condition $(a')$ is equivalent to the embedding 
$R_{+\cdot}^{\bar{b}}\subset 
\la(\bar{\si}\dot{w}_0^{\bar{b}})$.
The set $\la(\bar{\si}\dot{w}_0^{\bar{b}})$
does not contain $R_{+\cdot}^{\bar{b}}$ if and only if
at least one $\al_j\in R_{+\cdot}^{\bar{b}}$ 
is missing in the former set. Indeed, if all 
such $\al_j$ belong to this 
set then so do the roots that are
their positive linear combinations. 

If $\al_j\not\in \la(\bar{\si}\dot{w}_0^{\bar{b}})$, 
then $\al_j\in \la(w_0\bar{\si}\dot{w}_0^{\bar{b}})$.
Therefore $\al_j\in \varpi_b^\si$ due to
(\ref{varpibsi}) and because $(\al_j,\bar{b})=0$.
It gives the required.
\sq

Let us express the embeddings 
$$
R_{+\cdot}^{\bar{b}}\subset 
\la(\bar{\si}\dot{w}_0^{\bar{b}}),\ 
R_{+\cdot}^{b}\subset 
\la(\bar{\si}\dot{w}_0^{b}),
$$ 
using (\ref{barsidotw}). 
%\begin{align}\label{varpibsio}
%\varpi_b^\si = b (w_0\si \dot{w}_0^{b})^{-1}=
%%\bigl(w_0\dot{w}_0^{\varsigma({b})}\,
%(\si^\varsigma)^{-1}\bigr)\cdot(-\bar{b})
%\for \si^\varsigma=\varsigma(\si).
%\end{align}
Then ($ a'$) and ($ b'$) become 
equivalent correspondingly to 
\begin{align}\label{lengdbcnewa}
& (a''):\ \  
\si(R_{+\cdot}^{\,b})\,\widetilde{\cup}\,
\la(\si^{-1})
\,\supset\, \varsigma(R_{+\cdot}^{\bar{b}}) \and\\
\label{lengdbcnewb}
&(b''):\ \ 
\bar{\si}(R_{+\cdot}^{\,\bar{b}})\,\widetilde{\cup}\,
\la(\bar{\si}^{-1})\,
\supset \,\varsigma(R_{+\cdot}^{b}),
\end{align}
where $\,\widetilde{\cup}\,$ is the union where
the pairs $\{\tal,-\tal\}$ are removed.
 
These conditions can be simplified
if the following length conditions hold: 
\begin{align}\label{lengdotw}
&(\al):\  l(\dot{w}_0^{\bar{b}}\si^\varsigma)=
l(\dot{w}_0^{\bar{b}})+l(\si),\ \  
(\be):\  l(\dot{w}_0^{b}\bar{\si}^\varsigma)=
l(\dot{w}_0^{b})+l(\bar{\si}).
\end{align}
Then ($ a''$) and, correspondingly, 
($ b''$) become equivalent to 
\begin{align}\label{lengdbc}
& (a'''):\ \ \si(R_{+\cdot}^{\,b}) 
\supset \varsigma(R_{+\cdot}^{\bar{b}}), \ \ 
  (b'''):\ \ \bar{\si}(R_{+\cdot}^{\,\bar{b}})
\supset \varsigma(R_{+\cdot}^{b}).
\end{align}

For instance, let us check the equivalence
of ($ a''$) and ($ a'''$). Using that
$\la(\dot{w}_0^{\bar{b}})=R_{+\cdot}^{\bar{b}}$\,
and that $\si(\la(\si))=\-\la(\si^{-1})$,
$$
\la(\si^{-1})\cap \varsigma(R_{+\cdot}^{\bar{b}})=\emptyset
\ \Longleftrightarrow\  
(-\la(\si^\varsigma))\cap (\si^\varsigma)^{-1}
(R_{+\cdot}^{\bar{b}})=\emptyset
\Longleftrightarrow\  (\al).
$$

The requirements from (\ref{lengdotw})
simply mean that the left-hand side or
the right-hand side are reduced products
in the following transformation of
(\ref{barsidotw}):
\begin{align*}
\dot{w}_0^{\bar{b}}\varsigma(\bar{\si})\ =\ 
\si^{-1}\dot{w}_0^{\varsigma{b}}.
\end{align*}
This holds in many examples of minimal NGT,
though, generally, only for the left-hand side
or only for the right-hand side (corresponding to $(\al)$
or $(\be)$).

The setting from the next proposition
simplifies the
the construction significantly.  
%It does not cover all minimal NGT. For instance,
%the {\em parity corrections} for $\widetilde{C}_n$
%and $\widetilde{D}_n$ considered below require using
%more general $\si$. However, condition ($ a!b$) from
%the proposition is sufficient for $\widetilde{B}_n$.

\begin{proposition}\label{SIUPSB}
(i) Under the conditions $b\in P_+^{(k)}$ and  
$\bar{b}\in P_+^{(\bar{k})}$, let us assume that
\begin{align}
& (a!b):\ \ \si(\Ga^b_\cdot)= 
\varsigma(\Ga^{\bar{b}}_\cdot\,) \and 
\si(\al_k)=\varsigma(\al_{\bar{k}}).
\label{bisiprime}
\end{align}
Then  
\begin{align}\label{afactb}
&\bar{\si}\ =\ (\si^\varsigma)^{-1},\ \
\bar{{\si}}(\bar{b})\ =\ \varsigma(b),\ \ 
\dot{w}_0^{\bar{b}}\si^\varsigma\ =\ 
\si^\varsigma \dot{w}_0^{b^\varsigma},\\
&\la(\si^{-1})\,\cap\, \varsigma(R_{+\cdot}^{\bar{b}})
\ =\ \emptyset\ =\ 
\la(\bar{\si}^{-1})\,\cap\, \varsigma(R_{+\cdot}^{b}).
\notag
\end{align}
and the relations ($ a''',b'''$) from 
%(\ref{upsilonsharp})
(\ref{lengdbc}) are satisfied. 

(ii) Moreover,
\begin{align}\label{law0si}
&\la((\upsilon^\sigma_b)^{-1})=
\la(w_0\dot{w}_0^{b^\varsigma}(\si^\varsigma)^{-1})= 
R_+\setminus \bigl( R_{+\cdot}^{\bar{b}}\,\cup\, 
\la((\si^\varsigma)^{-1})\bigr),
\end{align}
where $R_{+\cdot}^{\bar{b}}\,\cap\,
\la((\si^\varsigma)^{-1}) =\emptyset$, 
i.e.,
the union is disjoint. Equivalently,
\begin{align}
l(\si)+l(\upsilon^\sigma_b)=l(\upsilon_b^\cdot)
\for \upsilon_b^\cdot=w_0\dot{w}_0^{b}\,,
\upsilon_b^\si=w_0\si\dot{w}_0^{b}\,.
\notag
\end{align}
The root  $\be=-w_0\si\dot{w}_0^b(\al_k)$ is
positive, equivalently,
$$
\al_k\in \la((\upsilon_b^\si)^{-1})=
\la(\dot{w}_0^b\si w_0).
$$
%and $p=(-b,\al_{k})$
\end{proposition}

{\em Proof.}
Relations (\ref{afactb}) obviously follow from
the assumptions from
(\ref{bisiprime}). Here 
$\la(\varsigma(\si)^{-1})\cap R_{+\cdot}^{\bar{b}}=
\emptyset$\ because $\varsigma(\si)^{-1}$ sends all simple 
roots from $R_{+\cdot}^{\bar{b}}$ to positive ones. 
Concerning  (\ref{law0si}),
\begin{align}\label{law0wb}
&\la(w_0\dot{w}_0^{b^{\varsigma}})\ =\ 
R_+\setminus \la(\dot{w}_0^{b^{\varsigma}})\ =\ 
R_+\setminus R_{+\cdot}^{b^{\varsigma}},\\
&\la(w_0\dot{w}_0^{b^{\varsigma}}(\si^\varsigma)^{-1})\ =\ 
R_+\setminus \bigl( R_{+\cdot}^{\bar{b}}\cup 
\la((\si^\varsigma)^{-1})\bigl).
\notag
\end{align}
It readily gives the desired length equality.

Let us check that $\be=
-w_0\si\dot{w}_0^b(\al_k)$ is positive. 
It suffices to assume that  
$(\be, \bar{b})=(-\al_k,-b)=(\al_k,b)<0$.
The positivity of $\be$ is
equivalent to the positivity of $\varsigma(\si(\al_k))=
\al_{\bar{k}}$.
Indeed, $\dot{w}_0^b(\al_k)=\al_k+c$, where
$c$ is a linear combination of the
simple roots from $\Ga^b_\cdot\,$. Therefore,
\begin{align*}
\be\ =\ -w_0\si\dot{w}_0^b(\al_k)\ =\ \varsigma(\si(\al_k))+c',
\end{align*}
where $c'$ is a linear combination of the
simple roots from  $\Ga^{\bar{b}}_\cdot\,$. 
Since $c'$ does not include
$\al_{\bar{k}}$, the root $\be$ is positive.
\sq
\medskip

If $b=0$ (the nonaffine case, which is allowed), 
then the condition ($ a!b$) implies that $\si$ must
be trivial. Therefore the $\si$\~extension becomes the 
{\em dot}\~extension in this case. Namely,
$\upsilon_{b=0}^\cdot= u_k=u_{\om_k}$ 
in the notations from (\ref{xwo}). 

In the case of trivial $\si=$\ id,
$$
\bar{b}=\varsigma(b),\ \al_{\bar{k}}=
\varsigma(\al_k) \and  (\varpi_{b}^\cdot)^{-1}=
\varpi_{\varsigma(b)}^\cdot \and
$$
$\be=-w_0\dot{w}_0^b(\al_k)=\varsigma(\al_k)+
\sum_{0\neq j\neq k} c_j\al_j^\varsigma$.

One can use the latter to analyze directly when
the corresponding 
$\varpi_b=b (w_0\dot{w}_0^{b})^{-1}$ is (represents)
a minimal NGT, i.e.,  when $\ga=\be+\al_{k}^\varsigma$ is a 
root. Generally, it leads to explicit conditions for the
coefficients $c_j$ for simple roots $\al_j$ neighboring $\al_k$
in the Dynkin diagram.

Concerning the positivity of $\be$ for minimal NGT
(which automatically results from the proposition),
given an affine minimal  NGT with $\be<0$ of types
$\widetilde{B},\widetilde{D}$,  one can apply the 
transposition of the mirrors and make $\be$ positive (see below).
Algebraically, this transformation corresponds to   
the symmetry $\al_0\leftrightarrow\al_n$. 
The case of $\widetilde{C}$ is analogous.
\medskip

\subsection{Universal construction}
Let us  address the root $\al_0$.
We use the notations and formulas from the previous section.

\begin{proposition}\label{BARBTH}
The following conditions
for $b,\si$ and $\bar{b},\bar{\si}$ are necessary and sufficient
for $\al_0\not\in \la(\varpi_b^\si)$
and $\al_0\not\in \la(\varpi_{\bar{b}}^{\bar{\si}})$
correspondingly:
\begin{align}\label{vthcond}
&(\bar{b},\vth)\le 0 \hbox {\  or \  }
(\bar{b},\vth)=1 \hbox{\ if\ \,} 
\dot{w}_0^{\bar{b}}(\vth)\not\in \la(\bar{\si}),
\notag\\
&(b,\vth)\le 0 \hbox {\ or\ \,}
(b,\vth)=1 \hbox{\ \,if\ \,} 
\dot{w}_0^{b}(\vth)\not\in \la(\si).
\end{align}
They imply that $b\not\in P_+$ and $\bar{b}\not\in P_+$
unless $b=0=\bar{b}$ or $b$ and $\bar{b}$ are (both) minuscule.
\end{proposition}
{\em Proof.} First of all, we note that  
$\dot{w}_0^{\bar{b}}(\vth)\not\in \la(\bar{\si})$
is equivalent to $\upsilon_b^\si(\vth)<0$, which
is sometimes easier to check.
Using the definitions from
(\ref{varpibsi}) and general formula (\ref{latutw}),
\begin{align}
&\varpi_b^\si\ =\ 
(\dot{w}_0^{b}\si^{-1} w_0)\cdot(-\bar{b})
%=(w_0\dot{w}_0^{b^\varsigma}(\si^\varsigma)^{-1})\cdot(-\bar{b}),
=(w_0\bar{\si}\dot{w}_0^{\bar{b}})\cdot(-\bar{b}),\notag\\
\label{lavarpib}
&\la(\varpi_b^\si)\ =\ 
\bar{b}\Bigl(
R_+\setminus \bigl(\dot{w}_0^{\bar{b}}(\la(\bar{\si})
\widetilde{\cup} R_{+\cdot}^{\bar {b}}\bigr)
\Bigr)\, \widetilde{\cup}\, \la(-\bar{b}),
%\bar{b}\Bigl(
%R_+\setminus \bigl(\si(R_{+\cdot}^{\varsigma(b)})
%\,\widetilde{\cup}\, \la((\si^\varsigma)^{-1})\bigr)
%\Bigr)\, \widetilde{\cup}\, \la(-\bar{b}),
\end{align}
where the modified union $\,\widetilde{\cup}\, $ includes,
by definition, removing all possible pairs $\{\tla,-\tla\}$. 

The root $\al_0=[-\vth,1]$ can appear in $\la(\varpi_b^\si)$
only from $\la(-\bar{b})$ because all other
roots there have positive nonaffine components.
It exists in $\la(-\bar{b})$ if and only if
$(-\bar{b},-\vth)=(\bar{b},\vth)\ge 1$; see (\ref{xlambi}). 
However, when it belongs to
$\la(-\bar{b})$, it can be still canceled by 
$[\vth,-1]=\bar{b}(\vth)=[\,\vth\,,\,-(\bar{b},\vth)\,]$ 
under the following conditions
$$
\vth\not\in \dot{w}_0^{\bar{b}}(\la(\bar{\si}))
\and (\bar{b},\vth)= 1.
$$ 
We use here that $\vth\not\in R_{+\cdot}^{\bar {b}}$
unless $b=0$. 

Now let us assume that $\bar{b}\in P_+$ under the
conditions for $\bar{b}$ from (\ref{vthcond});
the case of $b$ is analogous. Then
$(\vth,\bar{b})=(\vth,\varsigma(\bar{b}))>0$ 
unless $\bar{b}=0$. The conditions for $\bar{b}$
from (\ref{vthcond}) imply that $(\vth,\bar{b})=1$
and $\bar{b}$ is minuscule; so its
$\la$\~set is actually empty.
\sq
\smallskip 
%$$
%(\vth,\bar{b))=
%(\bar{\si}\dot{w}_0^{\bar{b}}(\vth),
%\bar{\si}\dot{w}_0^{\bar{b}}(\bar{b}))=
%(\bar{\si}\dot{w}_0^{\bar{b}}(\vth),
%\varsigma(b))
%$$

Note that if  ${\si}=$\ id, then $\bar{b}=\varsigma(b)$,
$(\bar{b},\vth)=(b,\vth)$ and
the conditions from (\ref{vthcond}) become
a single inequality $(b,\vth)\le 1.$  
The next proposition addresses the occurrence of $\al_0$ 
through the construction of
the elements $\pi_b$ from Proposition \ref{PIOM}.

\begin{proposition}\label{VTHPIB}
The elements $\hw\in \hW$ such that $\al_0$ is
a unique simple root in $\la(\hw)$ and, moreover, 
there exists only one simple root in 
$\la(\hw^{-1})$  are as follows:
\begin{align}\label{twpicr}
\hw=\pi_b \for
b=-m\om_i+\sum_{r\neq i}c_r\om_j,\ 
1\le i\le n, m\in \Z, c_r\in\Z_+. 
\end{align}
Then the endpoints of $\la(\pi_b)$
are unique. Namely, this sequence begins with  $\al_0$ and 
ends either with $\al_0$ for $m\le 0$ or
with $[-\al_i,\nu_i m]$ if $m>0$. Such elements $\hw$ never 
represent minimal NGT.
\end{proposition}   
{\em Proof.}
The element $\hw$ with $\la(b)$ such that $\al_0$ is its
unique simple root must be in the form $\pi_b$; it is
necessary and sufficient. 
Let us use (\ref{lapimin}):
$$
\la(\pi_b^{-1}) = \{ [\al,\nu_\al j]\in \tR_+,\ 
 -(b,\al^\vee)>j\ge 0 \}.
$$
We see that $(b,\al_i)$ can be strictly negative only for one
$i$ such that $1\le i\le n$. This gives the representation from
(\ref{twpicr}).

If $m\le 0$, then $2\al_0$ is not a root and the
corresponding $\pi_b$ is not a minimal NGT.
If $m>0$, then the roots  
$\vth, \al _i, \vth+\al_i$
cannot (all) belong to a root subsystem of $R$ of type $A_2$,
since $\vth$ is the maximal short root.  \sq 
\smallskip

The following theorem is essentially a combination
of the previous considerations.
We mainly focus on the general problem of finding 
adequate presentations for
$\hw\in \hW$ with non-movable endpoints in its 
$\la(\hw)$. As an application, it provides a convenient 
universal tool (i.e., for all root systems) for
managing the classification of the minimal NGT. 

Here and below we will constantly use that if 
$\hw$ is with non-movable endpoints or 
is a minimal NGT then so is $\hw^{-1}$.
Also,
\begin{align}\label{pitrans}
&\hbox {if\ }
\hw \hbox{\ is a minimal NGT, then so are\ \,}
\pi_r\hw \and \hw \pi_r,\  \pi_r\in \Pi.
\end{align} 
For instance, if $\hw=\pi_r\tw\in \hW$ represents
a minimal NGT, where $\tw\in \tW$, then
$\tw$ and $\pi_r(\tw)\pi_r^{-1}\in \tW$
are minimal NGT too.  The reduction $\hw\mapsto \tw$
and conjugations by $\pi_r$
can lead to quite non-trivial examples of 
minimal NGT (as words considered strictly inside
$\tW$) even if the initial $\hw$
is relatively simple.

\begin{theorem}\label{ALMOSTPOSITIVE}
(i)
Under the conditions from (\ref{lengdotwnew})
or, equivalently, (\ref{upsilonsharp}), 
the set $\la(\varpi^{\si}_b)$ 
contains only one simple nonaffine root, $\al_{\bar{k}}$,
the set  $\la(\bar{\varpi}^{\si}_b)$ 
contains only $\al_k$, where 
$$
\bar{\varpi}^{\si}_b\ =\ (\varpi^{\si}_b)^{-1}
\ =\ \varpi^{\bar{{\si}}}_{\bar{b}}.
$$
Conditions (\ref{vthcond}) guarantee
that $\la(\varpi_b^{\si})\not\ni 
\al_0\not\in \la(\bar{\varpi}_b^{\si})$.

(ii) Equivalently, under the same conditions,
the endpoints (the beginning and the end)
of the sequence $\la(\varpi_b^{\si})$ are 
non-movable, i.e., do not depend on the choice of the reduced 
decomposition of $\varpi_b^\si$. Namely, this sequence 
begins with $\al_{\bar{k}}$ and ends with
\begin{align}\label{tbevarpi}
&\tbe=-\varpi^{-1}(\al_{k})=[\be,p] \for 
\be\equal -\upsilon(\al_{k}),\ p=-(b,\al_k)\ge 0,
\end{align}
where we set $\upsilon=\upsilon^{\si}_b,\ 
\varpi=\varpi^{\si}_b$. It implies that
$b\not\in P_+ \not\ni \bar{b}$ unless $b=0=\bar{b}$;
thus, the dot-extension is actually needed only for
zero $b$.
%The root $\be$ is positive.

(iii) All elements $\hw$ with non-movable endpoints (both) 
are in the form  $\hw=\varpi_b^\si$ for 
$b\in P_+^{(k)}$ and $\bar{b}\in P_+^{(\bar{k})}$
subject to the assumptions from (\ref{lengdotwnew}) and
(\ref{vthcond}), and also provided that
$\la(\hw)\not\ni\al_0\not\in \la(\hw^{-1})$. The latter
condition always holds if we assume that 
$\ga=\al_{\bar{k}} +\be\in R$
for $\be$ from (ii).
In this case, $\tga\equal \al_{\bar{k}}+\tbe=[\ga,p]$ 
belongs to $\la(\varpi)$. All minimal NGT can be
obtained in this way.
\end{theorem}

{\em Proof.}
Part ($i$) is a combination of 
Proposition\ref{BARBTH}
and Proposition \ref{BARBLEN}. Formula 
(\ref{lavarpib}) reads now:
\begin{align}\label{lavarpibsi}
&\la(\varpi_b^\si)\ =\ 
\bar{b}\Bigl(
R_+\setminus \bigl(\dot{w}_0^{\bar{b}}(\la(\bar{\si}))
\cup R_{+\cdot}^{\bar {b}}\bigr)
\Bigr)\, \widetilde{\cup}\, \la(-\bar{b}),\notag\\
&\hbox{where\ }
\dot{w}_0^{\bar{b}}(\la(\bar{\si}))
\cap R_{+\cdot}^{\bar {b}}\ =\ \emptyset.
\end{align}
The only simple nonaffine
root from $\la(\varpi_b^\si)$ can be  
$\al_{\bar{k}}$
due to formula (\ref{upsilonsharp}). 
It comes from $\la(-\bar{b})$ when $p>0$. If
$p=0$, then this root will remain in this set only if
$\al_{\bar{k}}\not\in\dot{w}_0^{\bar{b}}(\la(\bar{\si}))
\cup R_{+\cdot}^{\bar {b}}$\,. This case can be managed
using a certain modification of (\ref{upsilonsharp}).
This can avoided because of the following argument.
At least one simple root must be present 
in $\la(\varpi_b^\si)$; therefore
it can be only $\al_{\bar{k}}$ since we excluded $\al_0$.
The case of $\la(\varpi_{\bar{b}}^{\bar{\si}})$
is entirely parallel.

Claim ($ii$) is a reformulation
of ($i$). Note that, generally, conditions 
from (\ref{lengdotwnew}),
(\ref{upsilonsharp}) and (\ref{vthcond}) make it
possible to reduce any claims about the last roots
(the ends) in the $\la$\~sequences under consideration
to the statements about the first roots (the beginnings). 
The interpretation of ($ii$) in terms of 
minimal NGT from ($ iii$) is obvious.

Concerning the implication $b\in P_+ \Longrightarrow 
b=0$, it was stated in Proposition \ref{BARBTH}. 
Recall that for dominant weights $b$, 
there is some flexibility with
picking $\al_{k}$ and $\al_{\bar{k}}$
orthogonal to $b$ and $\bar{b}$.
Without the reference to Proposition \ref{BARBTH}, 
the required implication 
follows immediately from the formula for $\tbe$.
Indeed, since $p$ is zero for $b\in P_+$, then $\tbe\in R_+$
and the triple $\{\tbe,\tga,\al_{\bar{k}}\}$ is nonaffine;
as such, it can be represented only by a nonaffine element.
The latter means that $b=0$.  

As for the completeness of our construction,
we need to check only the fact that
$\hw$ and $\hw^{-1}$ with non-movable ends can
be obtained from {\em almost dominant} weights
$b$ and $\bar{b}$. The other conditions we imposed
were actually necessary and sufficient.  

Let us assume that $\hw$ satisfies the
conditions from ($iii$) and set
$\hw=wc$ for $c\in P, w\in W$. If 
$(\al_k,c)> 0<(\al_{k'},c)$
for $k\neq k'$, then $\la(c)$ 
contains two simple roots $\al_k$ and
$\al_{k'}$. These roots cannot be canceled
in $\la(wc)=c^{-1}(\la(w))\widetilde{\cup}
\la(c)$ because all roots from 
$c^{-1}(\la(w))$ have positive nonaffine
components. This contradiction concludes
($iii$). \sq
\medskip

{\bf An example of type \mathversion{bold}{$D-E$}.}
Let us give a general construction applicable to 
all root systems
of types $D_{n} (n\ge 4), E_6, E_7, E_8$.
% which mainly
% addresses the condition $\al_{\bar{k}} +\be\in R$.
We take $\al_k=\al_4$ for $E_{6,7,8}$ and $\al_k=\al_{n-2}$
for $D_n$. The notation is from \cite{Bo};
this root has three neighbors in $\Ga$. Let $b=-\om_k$,
$\si=\ $id. For $D_n$, one has 
$b=-(\ep_1+\ep_2+\ldots +\ep_{n-2})\in P$.
Then the element $\bar{b}=\varsigma(b)$
coincides with $b$. Thus, $k=\bar{k}$ ($k=4$ for $D_4,E$).
Since $(b,\th)<0$, the conditions
from (\ref{vthcond}) are obviously satisfied.

One has: $\be=-w_0w_0^b(\al_k)=\sum_{j=1}^n {\al_j}$.
Here $w_0w_0^b$ sends 
$$ 
\ep_j\mapsto -\ep_{l-1-j} \for
j\le l-2,\ \ep_{n-1}\mapsto \ep_{n-1},\ 
\ep_{n}\mapsto -\ep_{n}.
$$ 
The relation $\al_k+\be\in R$ is not immediately clear
but readily follows from the tables of \cite{Bo}.
Thus, $\varpi=b w_0^bw_0$ represents the minimal NGT
$\{[1,\be],[1,\al_k+\be],\al_k\}$.

From the viewpoint of Theorem \ref{TYPECDAF} below,
the geometric classification theorem, this example
gives the simplest possible minimal NGT for $D_{n}$.
The corresponding configuration
has two bottom horizontal lines (the smallest possible
number); the other lines form a bunch of ``parallel lines"
with one reflection in the top mirror and $2$ in the bottom
one.  

If $n$ is odd here, then $\om_k\not\in Q$ and 
a {\em parity correction} is needed if we want to reduce
$\varpi$ to $\tW$. Namely, $\varpi'=\pi_1\varpi\in \tW$ 
and $\varpi''=\varpi\pi_1=(\varpi')^{-1}\in \tW$ 
are minimal NGT for $\pi_1=\pi_{\om_1}=\ep_1 s_{\ep_1}s_{\ep_n}.$
The corresponding weights $b$ in the decomposition 
$\tw=bw, b\in Q$ are
\begin{align}\label{d-example}
b'=2\ep_1-(\ep_2+\ldots +\ep_{n-2}),\ 
b''=-(\ep_1+\ldots +\ep_{n-3}+2\ep_{n-2}).
\end{align}
Note that the elements $\varpi_{b'}$ and $\varpi_{b''}$,
which are minimal NGT too, are different from
$\varpi'$ and $\varpi''$, although the corresponding
weights coincide.
\sq
 
%For $D_{2l+1}$, one can take here
%(among other choices) the weight
%$b=2\ep_1-(\ep_2+\ldots+\ep_{n-2})$. It provides an
%example that is beyond the following corollary.
%See Figure \ref{cngt2}.
\smallskip

The theorem reduces the classification of
all elements $\hw$ with non-movable endpoints
of their $\la$\~sequences to {\em finitely 
many} verifications. The analysis is involved
for the exceptional root systems. 
However, we expect that the classification
of minimal NGT (a subclass of all $\hw$ with
non-movable $\la$\~endpoints) 
is not that ramified. For instance,
the description of all minimal NGT satisfying Proposition 
\ref{SIUPSB} seems quite doable (although this setting
is not sufficient for all of them). It is similar to the
verifications in the example above. 
\smallskip

Combining this
proposition for the simplest $\si=\ $id with the
natural extensions to greater root systems we come to the
following corollary, essentially, sufficient to obtain
all classical affine minimal NGT.  
\smallskip

\begin{corollary}\label{SUBSYSNGT}
Let $R'$ be a root subsystem of $R$ such that the
corresponding Dynkin diagram $\Ga'\subset \Ga$ is
connected and contains $B_3,C_3$ or $D_4$ as a subdiagram, 
$$
\varpi_{b'}^\cdot=b(\upsilon_{b'}^\cdot)^{-1}\in
\hW' \for \upsilon_{b'}^\cdot=w_0'\dot{w}_0^{b'},
$$
subject to the conditions for $b'\in P'$
from the theorem (with $\si'=\,$id). We 
linearly extend the embedding $Q'\subset Q\subset \R^n$ 
to $P'\subset \R^n$. 

(i) Let us assume that $b'\in P$ and that there
exists $\,b\in P\,$ such that $b-b'$ is a linear
combination of the fundamental weights $\om_j$ for
$\al_j\not\in \Ga'$, provided the conditions
\begin{align}\label{bprimext}
&b\in P_+^{(k')},\ \, \bar{b}\equal
-\upsilon_{b'}^\cdot(b)\,\in\, P_+^{(\bar{k}')},\\ 
&(b,\vth)\,\le\, 0\, \ge\, (\bar{b},\vth)\, \for \,
\vth\in R_+,\notag  
\end{align}
where $k', \bar{k}'$ correspond to $b', \bar{b'}$.
Then $\hw\equal b(\upsilon_{b'}^\cdot(b))^{-1}$
represents a minimal NGT for $\widetilde{R}$.
%the endpoints of its $\la$\~sequence 
%remain unchanged and non-movable for $\tW$ and $\widetilde{R}$

(ii) If $b=b'$ here (the extension by zero),
then $\bar{b}=-\varsigma'(b)$
for $\varsigma'$ defined for $R'_+$. In this case
one must check that 
$(\varsigma'(b'),\al_m)\ge 0$ 
for any simple root $\al_m$ neighboring
$\Ga'$ in $\Ga$ (two may occur for $E_{6,7,8}$).
Also, if
$\al_0=[-\vth,1]\in \tR$ is connected inside the
affine diagram $\widetilde{\Ga}$ for $\widetilde{R}$
with one of the vertices of $\Ga'$ by a link, then
the conditions  
$(b,\vth)\le 0\ge (\varsigma'(b),\vth)$ must hold.
Such $\hw$ represents a minimal NGT. 
\end{corollary}
{\it Proof}. 
Let us  demonstrate 
that the endpoints of $\hw$ coincide with those of
$\varpi_{b'}^\cdot$ (and are unique).
Representing $\hw=
(\upsilon_{b'}^\cdot)^{-1}(-\bar{b})$,
we see that new {\em simple} nonaffine roots can appear only 
due to the set $\la(-\bar{b})\subset \widetilde{R}$.
If $b=b'$ then  $\bar{b}=-\varsigma'(b)$ and 
they have to be among the neighbors $\al_m$
of $\Ga'$ in $\Ga$ such that $(\al_m,-\varsigma'(b'))>0$.
The analysis of $\al_0\in \widetilde{R}_+$ is 
straightforward using (\ref{vthcond}). \sq
\smallskip

We note that the special case $b=b'$ is actually covered
by Proposition \ref{SIUPSB}.
%The {\em parity corrections} considered below require
%more general extensions $b\neq b'$ (and even this
%modification is not sufficient to cover them all).
The {\em dots} in this corollary can be ignored
for $b'\neq 0$. If $b'=0$, i.e.,
the initial $\varpi_{b'}$ is non-affine, then
$b'$ can be extended only by zero due to the
relations in terms of $\vth$.

One can combine the transformations $\hw\mapsto \pi_r\hw$
and $\hw\mapsto \hw\pi_r$ for $\pi_r\in \Pi$ with the
construction from the corollary;
see (\ref{pitrans}).
Generally, the resulting elements 
will be ``new", i.e., not covered 
this corollary for any proper $b'$ and their extensions
$b$. Multiplication by $\pi_r$ here may change the 
centralizer of $b$ and result in significant changes of 
the $\sigma$\~elements of the theorem.

Furthermore, using the above transformations
and automorphisms of $\tGa$, one obtains 
all affine minimal NGT of classical types.

\comment{
The justification of the following corollary is
very much similar to this one; we will omit it.
It allows us to increase $R^b_k$ under sufficiently
relaxed conditions, and includes the {\em parity
corrections} for the classical root systems 
considered below.

\begin{corollary}\label{RINRPRIME}
Let us assume that $R^b\subset R'\supset (R^{\bar{b}})^\varsigma$ 
for a root subsystem $R'\subset R$ such that
\begin{align}\label{RinRpr}
&\bar{b}=\varsigma(w_0'(b))\in P_+^{(\bar{k})} \and
(\bar{b},\vth)\le 0\ge (b,\vth)
\end{align}
for the maximal element $w_0'$ associated with $R_+'$.
Then $\si=w'_0w_0^b$ and the sequence
$\la(\varpi)$ for  $\varpi=\varpi_b^{\si}=
bw_0'w_0$ has non-movable endpoints
$\al_k$ and $[\be,p]$, where $\be=\varsigma(w_0'(\al_k)),$\,
$p=-(b,\al_k)>0$. Provided that $\ga=\al_k+\be\in R$, it
represents a minimal NGT.\sq
\end{corollary}
\medskip
}

\smallskip
The next corollary is about applications of our construction
to {\em nonaffine} minimal NGT.
\smallskip

\begin{corollary}\label{NONAFFNGT}
Let $\varpi$ be the element from part (ii) of 
Theorem \ref{ALMOSTPOSITIVE} satisfying the
assumptions there. We require the positivity
of $\be$; for instance, the setting of
Proposition \ref{SIUPSB} is sufficient.
Using the decomposition
$\varpi=\pi_{\varpi}u_{\varpi}$ from Proposition \ref{PIOM},
the element $u_{\varpi}\in W$ represents a nonaffine minimal 
NGT $\{\be,\ga,\al_{\bar{k}}\}$ under the following condition. 
For any end $\be'$ of the sequence 
$\la(u_{\varpi})$, the root $[\be',-(\be',\bar{b})]$ must be 
from $\la(\varpi)$, equivalently,
$(\si^\varsigma)^{-1}(\be')>0$.
All classical nonaffine minimal NGT can be obtained 
as  $u_{\varpi}$ for $b=0$ and $\si=\,$id, possibly, with
further embedding into a greater root system via
Corollary  \ref{SUBSYSNGT}.
\end{corollary}
{\em Proof.}
The claim concerning
the classical root systems can be checked by inspection
(see \cite{CS} and below). Concerning the general
statement, let us begin with clarifying the structure 
of $\la(\varpi)$. If $\si$ is known (trivial or relatively simple),
then  the calculation of $\la(u_{\varpi})$ for 
$\varpi=\varpi_b^\si$ becomes sufficiently explicit.
Let us use (\ref{lavarpib}) in the following form:
\begin{align}\label{lavarpibb}
&\la(\varpi_b^\si)\ =\ 
\bar{b}\Bigl(
R_+\setminus \bigl(\dot{w}_0^{\bar{b}}(\la(\bar{\si})
\cup R_{+\cdot}^{\bar {b}}\bigr)
\Bigr)\, \widetilde{\cup}\, \la(-\bar{b}).
\end{align}

Here the nonaffine roots can appear due to $\la(-\bar{b})$
and because there can exist $\al\in R_+$ satisfying
$(\al,\bar{b})=0$ not from $R^{\bar{b}}_{+\cdot}$.
Thus,
\begin{align} \label{lapivarpi}
\la(u_{\varpi})=
\{\al\in
R_+\setminus 
&\bigl(\dot{w}_0^{\bar{b}}(\la(\bar{\si})
\cup R_{+\cdot}^{\bar {b}}\bigr)
\,\mid\, (\al,\bar{b})=0\}\notag\\
&\widetilde{\cup}\,
\{\al\in R_+\,\mid\, (\al,\bar{b})<0\},
\end{align}
where the set in the second line
obviously does not intersect the first one.

Apart from the roots with $(\al,\bar{b})=0$,
the description of $\la(\varpi)$ and
$\la(u_{\varpi})$ becomes simpler. Let us use
the relation (\ref{barsidotw})
\begin{align}\label{barsidotww}
&\bar{\si}\dot{w}_0^{\bar{b}}\ =\ 
\dot{w}_0^{b^\varsigma} (\si^\varsigma)^{-1}.
\end{align}

First, we describe the 
set of all $[\al,\nu_\al j]$ in $\la(\varpi)$
with $\al>0$ such that $(\al,\bar{b})=-q<0$. 
The inequality for $j$ is $0\le\nu_\al j\le q$
if $\al$ is {\em not} from $\la((\si^\varsigma)^{-1})$. 
We use that $[\al,q]$ belongs to the first set 
in the union
from (\ref{lavarpibb}). Otherwise, i.e.,
when  $\al\in \la((\si^\varsigma)^{-1})$,
the inequality is $0\le\nu_\al j<q$.
Summarizing,
\begin{align}\label{lavarpos}
\la(\varpi)&\,\cap\, [\al,\Z_+]\ = \\
&\left\{\begin{array}{ccc}
\{[\al,\nu_\al j]\,\mid\, 0\le\nu_\al j\le q \}
\hbox{\ \,if\ } \al\not\in \la((\si^\varsigma)^{-1}),&&\\
\{[\al,\nu_\al j]\,\mid\, 0\le \nu_\al j< q \} 
\hbox{\ \,if\ } \al\in \la((\si^\varsigma)^{-1}).&&
\end{array}\right.\notag
\end{align}

Second, let $\al>0$ and $(\al,\bar{b})=q>0$.
In contrast to the previous case, 
$\si$ may lead to diminishing the corresponding
part of $\la(\varpi)$:
\begin{align}\label{lavarneg}
\la(\varpi)&\,\cap\, [-\al,\Z_+]\ = \\
&\left\{\begin{array}{ccc}
\{[-\al,\nu_\al j]\,\mid\, 0<\nu_\al j< q \}
\hbox{\ \,if\ } \al\not\in \la((\si^\varsigma)^{-1}),&&\\
\{[-\al,\nu_\al j]\,\mid\, 0<\nu_\al j\le q \} 
\hbox{\ \,if\ } \al\in \la((\si^\varsigma)^{-1}).&&
\end{array}\right.\notag
\end{align}
\smallskip

Let us now verify the claim concerning $\la(u_{\varpi})$.
Recall that $\tbe=-w_0\si\dot{w}_0^b(\al_k)=[\be,p]$\, for
$p=-(b,\al_k)=-(\bar{b},\be).$ Thus the relation we imposed
on $\be'$ is satisfied for $\be$.
One can assume that
$p>0$, which excludes the case of zero $b,\bar{b}$.

In the reduced
decomposition $\varpi=\pi_{\varpi}u_{\varpi}$ for
$\varpi=\varpi_b^\sigma$, the set $\la(\pi_{\varpi})$
does not contain nonaffine roots, so the first
(and simple) nonaffine root in the $\la$\~sequence
of $\varpi$ must be from $\la(u_{\varpi})$.
Thus, the condition that the beginning of the
$\la(\varpi)$ is not movable (under the Coxeter
transformations) is {\em equivalent} to the corresponding
property of $u_{\varpi}$ provided that 
$\al_0\not\in \la(\varpi)$.

Using the positivity of $\be$ from $\tbe$, we conclude
that $\be\in \la(u_{\varpi})$. So does $\ga=\be+
\al_{\bar{k}}$. By assumption,
if $\be'\neq \be$ is an end of the sequence $\la(u_{\varpi})$,
then it can be ``lifted" to the root  $\tbe'=[\be',p']$
from $\la(\varpi)$, where $p'=-(\bar{b},\be')$. 

We will use
the following interpretation of the (left) {\em ends} $\tbe'$
of a $\la$\~sequence (see \cite{C0}). There must be no 
decompositions $\tga'=\tbe'+\tal'$ in terms of
positive roots (including the pure imaginary ones) where
$\tga'$ belongs to a given $\la$\~sequence and 
$\tal'$ does not. Also, there must be no decompositions
$\tbe'=\tal'+\tal''$ with $\tal'$ and $\tal''$ from this
$\la$\~sequence. These conditions are necessary and sufficient.

The decomposition 
$\tbe'=[\be',p']=[\al',r']+[\al'',r'']$ in terms of the roots
from $\la(\varpi)$ with positive $\al'$ and $\al''$
is impossible, since $\be=\al'+\al''$ contradicts
the assumption that $\be$ is a (left)
end of $\la(u_{\varpi})$. However, it 
can be in the form $\tbe'=[\al',r']+[-\al'',r'']$ 
for positive $\al'$ ($r'\ge 0$) and 
$\al''$ ($r''>0$). 
Then $\al'\in \la(u_{\varpi})$
and $\al''\not\in \la(u_{\varpi})$. Hence,
$\al'=\be+\al''$ and $\be$ cannot be an end of 
$\la(u_{\varpi})$ since it must appear before $\al'$.

Next, let us consider the case when 
$\la(\varpi)\ni [\ga',q']=
\tbe'+[\al',r']$  for $[\al',r']\not\in \la(\varpi)$.
If $\al'>0$ then $\be'$ cannot be the end.
If $\al'=-\al$ for $\al>0$ and $\ga'>0$, then
$\ga'\in \la(u_{\varpi})$ and  $\be=\ga'+\al$.
Therefore, $\al\not\in \la(u_{\varpi})$.
One can assume here (and below) that $r'=1$. 
The relation $\al\not\in \la(u_{\varpi})\not\ni [-\al,1]$
results in $(\bar{b},\al)=0$. Therefore, $(\ga,-\bar{b})=
(\be,-\bar{b})=p'$. It contradicts to
$q'=p'+1>p'$. See (\ref{lavarpos}); note that this
formula gives the reformulation of the 
assumption concerning the ends $\be'$ of $\la(u_{\varpi})$ 
we imposed.
 
The remaining case is when $\al'=-\al<0, r'=1$ and
$\ga'=-\ga<0$. Then $\tbe'+[\ga,-1-q']=[\al,-1]$.
If $\al \in \la(u_{\varpi})$, the root $\be$ must
appear in this sequence before $\al$, which is impossible.
Therefore, $(\bar{b},\al)=0$ and $[-\ga, -(\ga,-\bar{b})+1]$
belongs to $\la(\varpi)$; a contradiction.
See (\ref{lavarneg}) for detail.
\sq
\smallskip

It is of interest to explore the procedure from
Corollary  \ref{NONAFFNGT}
for obtaining nonaffine minimal NGT for the 
exceptional root systems; we have no claims so 
far concerning the completeness.
\smallskip

\setcounter{equation}{0}
\section{NGT of type
\texorpdfstring{\mathversion{bold}$B$}{\em B}}
The root system $\widetilde{B}_n (n\ge 3)$ is the key.
Due to our choice of $\vth$ (it is the maximal {\em short}
root; the twisted case), the corresponding 
affine Dynkin graph $\widetilde{\Ga}$,
$\Ga$ extended by $\al_0=[1,-\ep_1]$,
is the one from the $C$\~table of \cite{Bo} where
all the arrows are reversed. Concerning the 
normalization of the inner product,
our one is $(\ep_i,\ep_j)=2\de_{i,j}$ for the Kronecker
delta in terms of the basis $\{\ep_j\}$ from \cite{Bo},
 
The lattice $Q$ is generated by $\{\ep_j\}$. 
The action of $\tW$ in $\R^{n+1}=[z,\ze]$ and
$\R^n\in z=(z_1,\ldots,z_n)$ (see (\ref{ondthr}) 
and (\ref{standact})for $\xi=1$)
is as follows:
\begin{align*}
&s_0([z,\ze])\ =\ [(-z_1,z_2,\ldots,z_n),\ze-2z_1],\\  
&\,\ep_j(\!(z)\!)\ =\ z+\ep_j \for z=(z_1,\ldots,z_n).
\end{align*}

We will use the involution 
of $\widetilde{\Ga}$ transposing $\al_0=[-\ep_1,1]$
and $\al_n=\ep_n$; it will be denoted by
$\imath_B$. It coincides with the conjugation by 
$\pi_{n}\in \Pi$. 
\smallskip

For $\widetilde{C}_n$,
the lattice $Q$ is generated by $\ep_i\pm \ep_j$
including $2\ep_i$. Also, $\al_0=[-\vth,1]$
for $\vth=\ep_1+\ep_2$ and $(\ep_i,\ep_j)=\de_{ij}$. 
Accordingly, 
\begin{align*}
&s_0([z,\ze])=[(-z_2,-z_1,z_3,\ldots,z_n),\ze-z_1-z_2],\\  
&(\ep_i\pm\ep_j)(\!(z)\!)=z+\ep_i\pm\ep_j 
\for z=(z_1,\ldots,z_n).
\end{align*}

The affine
Dynkin diagram for $\widetilde{C}_n$ is
the one from \cite{Bo} for $B$ with all arrows reversed.
Its involution, transposing 
$\al_0=[-\ep_1-\ep_2,1]$ and $\al_1=\ep_1-\ep_2$
and fixing the other simple roots, 
will be denoted by $\imath_C$. In terms of $\ep_i$,
it sends $\ep_1\mapsto -\ep_1$ and leaves all other
$\ep_j$ unchanged. 

The lattice $P$ for $C$ coincides with $Q$ for $B$.
Sometimes we will denote $\, \al_0,\vth,s_0\,$ 
of type  $C$
by $\,\al_0',\vth',s_0'\,$ (the same for the related objects)
to avoid confusions with those defined for $B$.
For instance, the element $s_0$ from $B$ can be
treated as an element from $\hW=W\lsmash P'$  
defined for $C$. Namely, $\Pi'=\{\hbox{id},s_0\}$, i.e.,
$\pi_1'=s_0$. We see that $s_0$ induces $\imath_C$.  
\medskip

\subsection{Configurations of type 
\texorpdfstring{\mathversion{bold}$B$}{\em B}}
\begin{figure}[htbp]
\begin{center}
\vskip 1.5in
\hskip -.5in
\includegraphics[scale=0.45]{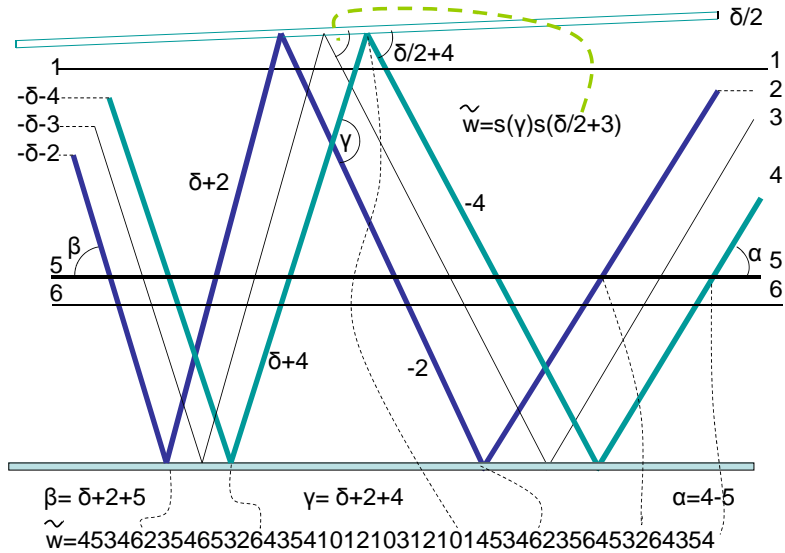}
\vskip -2.in
\caption{A basic affine NGT of type B}
\label{bngt}
\end{center}
\end{figure}

Let us begin with a simple typical example of 
minimal affine NGT of type $B$
presented in Figure \ref{bngt}.

There are $n=6$ lines there which intersect and 
also experience reflections in the two {\em mirrors}. The
bottom one will be always made parallel to the $x$\~axis,
the top one makes the angle $\de/2$ with this axis.

Here and further by a {\em line} we mean a 
piecewise linear zigzag line which is the result 
of reflections of the initial line in 
the {\em mirrors}. The latter will be referred to as 
the {\em bottom nonaffine} mirror and the {\em top affine} 
one.

Almost always we consider only the portion of such zigzag
lines trapped between the vertical lines at the
beginning and at the end of the graph. 
 
{\em Configurations} are defined as sets of (zigzag) lines
between a given pair of vertical lines (the beginning and the
end) where the
triple intersections and double reflections are not allowed.

The {\em initial} angles the lines make with the $x$\~axis 
(counterclockwise) will be denoted
by $\ep_1,\ep_2,\ldots,\ep_n$; we simply use $1,2,\ldots,n$
in the graphs.
They are numbered from top to bottom;
accordingly, the initial angle between
line $i$ and line $j$ for $j>i$ is denoted $i-j$.
We read the configuration from right to left, so its
beginning is the extreme right vertical line.

Let us list and interpret the geometric features of
configurations aiming at establishing connections 
with the algebraic theory 
of $\widetilde{B}_n$ and the corresponding $\tW$. 
\medskip

\vskip 0.5in
\medskip

{\bf Geometric features of configurations.}
%\smallskip

($ i$) The construction of 
$\tw\in \tW$ from a given configuration 
is explained in the figure;
see also \cite{CS} and \cite{Ch5}. More formally, 
we intersect the (zigzag) lines with the extreme
left vertical line and read the intersection
points from top to bottom, forming the sequence
of the {\em absolute} 
angles, which the lines make with the $x$\~axis. This sequence
can be uniquely represented as follows:
\begin{align}\label{debw}
&\de(b_1,b_2,\ldots,b_n)+w(1,2,\ldots,n)
\hbox{\ for proper }  b_i\in \Z \and 
w\in W.
\end{align}

Then the element
$\tw\in \tW$ (of type $\widetilde{B}$) associated with
the configuration is  defined as
the product $\tw=bw$. Here the vector
$(b_1,\ldots,b_2)$, which equals  $(0,-1,-1,-1,0,0)$ in the 
figure, is naturally identified with the weight
$b=\sum_{i=1}^n b_i\ep_i\in Q$.
Recall, that in the nonaffine theory of classical
Weyl group of types $B,C,D$ elements $w\in W$ 
are naturally 
identified with {\em permutations with signs}. 
For this particular configuration, $w=(1,-4,-3,-2,5,6)$.

Notice that the ``unit" here is $\de$
(not $\de/2$ as in the interpretation of the affine roots);
only integral multiples of $\de$ appear in the angles.
For instance, the vector of the absolute angles
after the event $s_0$ is $\de(1,0,0,\cdots)+(-1,2,3,\cdots)$. 
Thus the corresponding $b$ equals $\ep_1=\vth, w=s_{\vth}$, 
which matches the formula $s_0=\vth s_{\vth}$. 

As an exercise, check that $\tw$ from the figure
can be represented as a product of two pairwise 
commutative reflections $s_{[1,\ep_2+\ep_4]}$ and 
$s_{[1,\ep_3]}$.
\smallskip
    
($ ii$) The sequence of projections of the intersection
points and the reflection points onto
the $x$\~axis gives the {\em reduced decomposition}
of $\tw$ corresponding to a given configuration.
{\em We always assume that these projections are distinct}.  
Then their number equals the length $l(\tw)$.
The simple reflections $s_i$ ($0\le i\le n$)
associated with the corresponding {\em simple events},
the intersections and the reflections,
are determined on the basis of the {\em local line
numbers} (always counted from top to bottom) at the moment
of the event.

For disconnected events (corresponding to pairwise
commutative $s_i$ and $s_j$) we can of course change
the order of the projections arbitrarily; we 
do it constantly in the figures.

Note that if ``pseudo-lines" are allowed here,
then all reduced decompositions of a given
$\hw$ can be obtained in this way. Pseudo-lines are 
essentially the curves with 
one-to-one projections onto the $x$\~axis that are
allowed to intersect no greater than one time if no 
reflections are involved.
\smallskip

($ iii$)
Next, the angles $\al+(\de/2)j$ between the lines 
will be treated as the affine roots 
$[\al,j]\in \tR$ (type $\widetilde{B}$).
 The angle is always calculated 
{\em counterclockwise} and
{\em before the event}, i.e., as the difference of the 
absolute incoming angles, the upper one minus the lower one.
The events are intersections or reflections.
The angles with the mirrors are taken for the
reflections, namely, the absolute
angles of the mirror are $\de/2$ for the top one and 
$0$ for the bottom one. 
%\smallskip

The angles correspond to {\em positive} affine
roots, for instance, $\de/2$ always occurs with 
a non-negative coefficient (even for the intersections). 
The collections of the corresponding
angles considered from right to left 
constitute the $\la$\~sequences $\la(\tw)$
of a given reduced decomposition of $\tw$.
If {\em pseudo-lines} are allowed instead
of (straight) lines we consider, then all 
$\la$\~sequences can be obtained in this way. 
 
\smallskip
($ iv$) The action of $\tw$ on the angles is {\em dual} to
the affine action from (\ref{debw}). Practically, the
image of $\ep_i$ considered as a root is the resulting
angle of this line where index $i$ 
is replaced by the {\em local} number
of this line after the event (counted from top to bottom).
 
For instance, $\tw(\ep_2-\ep_5)$ in the figure
under consideration is 
$\tw(2-5)=(-\de-4)-5=-\de-4-5$ treated as 
the affine root $[-1, -\ep_4-\ep_5]$. It is negative,
so $2-5$ belongs to the list of the angles of this
configuration.
%\smallskip

Notice that the action of the lattice $P$ 
(of type $\widetilde{B}$) 
requires an extension of the basic events by $\pi_n$ 
transposing the affine Dynkin diagram $\widetilde{\Ga}$.
Recall that $\pi_n$ is the only non-trivial element of $\Pi$.
This event has no angle and does not contribute to the
$\la$\~sequences, although it of course transposes the
line numbers and influences the angles afterwards.   

Geometrically, let us assume that the mirrors are two
generatrix lines of a circular $2$\~dimensional cone;
then the configurations under consideration will
belong to the one of the two halves of this cone. 
The reflection in the middle line between the mirrors in the 
{\em other half of the cone} naturally represents $\pi_n$.   
It transposes the mirrors and the corresponding
lines between them; we denote it by $\imath_B$.
%\medskip

\subsection{ 
\texorpdfstring{\mathversion{bold}$B$}{\em B}-positive
minimal NGT} \label{sec:B-NGT}
We need to introduce some terminology.
%\smallskip

A collection of {\em neighboring parallel lines} 
will be called a {\em bunch of lines}.
The lines from a bunch are obtained from
each other by (piecewise) parallel translations 
(adjusted to the mirrors). 

Actually, by {\em parallel},
we mean here and below
{\em combinatorially parallel}, i.e., the
lines that ``behave" as parallel and may intersect only due 
to the reflections (within the range where they are considered), 
{\em We always assume that any
bunch is maximal possible} in a given configuration.

The lines from one bunch have the same 
numbers of top and bottom reflections. 
By horizontal, we mean the lines that are {\em 
parallel} ({\em combinatorially parallel}, to be exact)
to the corresponding mirrors; then these numbers
are zero. The {\em $t$\~number} of a line is defined 
as the number of top reflections; 
  
A natural generalization of the minimal NGT from 
Figure \ref{bngt}
is given in terms of the following data: 
\smallskip

($ a$) the integers $u\ge 0, v\ge 1$ such that $m\equal
n-u-v\ge 2$, which are the numbers of top and bottom
horizontal parallel lines neighboring (the right ends and 
the left ends) the corresponding mirror; 

($ b$) a decomposition $m=p_1+p_2+\ldots+p_r$ for positive 
integers $p_j$ such that  $p_{r}\ge 2$, which give the
numbers of lines in the consecutive non-horizontal 
bunches (counted from top to 
bottom with respect to the right ends); 

($ c$) a sequence of non-negative integers 
$0\le t_1<t_2<\ldots<t_r$, which are the $t$\~numbers of
the corresponding non-horizontal bunches;

($ d$) also, the number of the bottom reflections
is assumed $t+1$ for the bunches and $t$\~numbers in ($ b$); 
\smallskip

The data from ($a,b,c)$ determine the configuration uniquely
due to assumption ($d$).

\begin{figure}[htbp1]
\begin{center}
\vskip 1.5in
\hskip -0.5in
\includegraphics[scale=0.45]{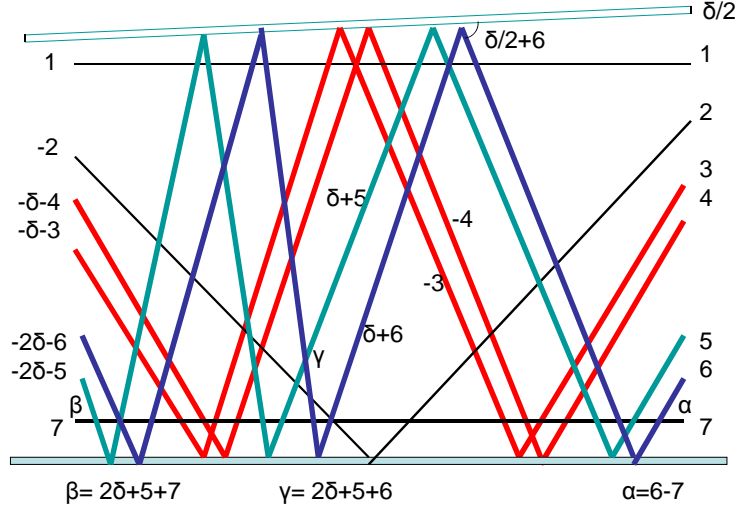}
\vskip -2.in
\caption{Basic affine NGT of type B}
\label{bngt1}
\end{center}
\end{figure}

Geometrically, the horizontal bunches can be plotted arbitrarily 
close to the corresponding mirrors; the lines in one bunch
can be plotted arbitrarily close to each other.
In Figure \ref{bngt1}, there are $1+1$ {\em horizontal}
bunches near the top mirror
and the bottom mirror 
(each with one line), namely, $\{1\}$ and $\{7\}$;
then $t_1=0,t_2=1,t_3=2$ 
for the bunches $\{2\},\{3,4\},\{5,6\}$.

%\vskip 1.5in
\medskip

The $t$\~number can be zero in our construction not only
for the {\em horizontal} bunches. 
The first bunch of lines from ($ b$) is allowed to 
have $t_1=0$.  
The presence of at least one {\em horizontal} bottom 
line ($v>0$) is required. Also, 
the second bunch counted from the bottom,
i.e., the first bunch from ($ b$),
must contain at least two lines
($t_r\ge 2$). 

The first and the last lines
from this bunch and the highest line
in the bottom {\em horizontal} bunch 
will be exactly those responsible for 
producing the minimal NGT $\{\tbe,\tga,\tal\}$
in the theorem below.

Note that the $\tw$\~elements corresponding to
different configurations under consideration
may have coinciding weights $b$. It occurs if and only 
if they have the same total number
of lines with $t=0$ due to a redistribution of lines between the
top horizontal branch and the one with $t_1=0$.
\smallskip

This construction will be referred to as
the {\em $B$\~positive construction}; accordingly,
such minimal NGT will be called {\em $B$\~positive}.
This name reflects the fact that the
nonaffine component $\be$ of the root $\tbe$ is always
positive in this construction. All minimal
NGT with positive $\be$ can be obtained in this way; we
come to the following theorem.

\begin{theorem}\label{TYPEBAFF}
(i) Any minimal affine NGT $\tw\in \tW$
for the (twisted) root system $\widetilde{B}_n (n\ge 3)$ 
is either given by the $B$\~positive construction in terms of  
($ a,b,c,d$) or can be obtain from a $B$\~positive minimal NGT 
by applying the automorphism $\imath_B$ (transposing
the top and the bottom mirrors). 
All such $\tw$ are involutive.

(ii) 
The $B$\~positive $\tw$ are covered by the construction of
Corollary \ref{SUBSYSNGT}, (ii) for $k=\bar{k}=v+1$,
where $v+1$ is the greatest line number in
the bunch of lines for the last $t_r$ 
(second from the bottom). The graph $\Ga'$ is obtained
from $\Ga$ by removing the vertices $\al_1,\ldots,
\al_u$, geometrically, by removing the top horizontal bunch 
(if present).
\sq
\end{theorem}
\medskip

\subsection{Proof}
We consider the configurations of the lines $L_i$
discussed above and representing elements the $\tw\in \tW$ for
$\tW$ of type $\widetilde{B}$. The lines are numbered
at the beginning (for the extreme right 
value of $x$). Each of $L_i$ is characterized by 
the number of top reflections $t_i$ and the number
of bottom reflections $b_i$. More exactly, this numbers
determines the type of $L_i$ uniquely if $t_i\neq b_i$;
otherwise, one needs to know which reflection
(the top or the bottom one) occurs the first.

Let us begin with the following general observation.

\begin{lemma}\label{BTINEQ}
Let the lines $L_i$ and $L_{i+1}$ be neighboring 
in the configuration corresponding to
$\tw\in \tW$ of type $\widetilde{B}$. 

(a) If the first reflection of line $L_{i}$ 
is in the bottom mirror and either $i=n$, or
$b_i<b_{i+1}$ or $t_i<t_{i+1}$ for $i<n$,
then the element $s_i$ can be made the 
beginning of the reduced decomposition of $\tw$. 

(b) Similarly, 
$s_i$ can be made the beginning of the reduced decomposition 
of $\tw$ if $L_{i+1}$ begins with the top 
reflection and either $i=0$ or, in the case of 
$i>0$, $t_{i}>t_{i+1}$ or $b_{i}>b_{i+1}$.
\end{lemma}
{\em Proof.} It suffices to check ($a$); also, the case
$i=n$ is obvious. The geometric assumptions from ($a$)
ensure that the angle $\ep_i-\ep_{i+1}$ occurs somewhere
in such configuration.
Indeed, the first reflection of line $L_{i+1}$ (if any) 
can be only in the bottom mirror. Then lines $L-i$ and
$L_{i+1}$ can be made 
``parallel" (i.e., with the intersections only due to 
their reflections) until the first intersection. Since the lines
have experienced the same number of the bottom and top reflections
before the intersection, the angle between them has to be 
$\ep_i-\ep_{i+1}$.  
This angle corresponds to the simple root $\al_i$;
therefore it can be made the first upon a proper transformation 
of the configuration.
\sq

\begin{lemma} \label{REDTO3}
The statement of Theorem \ref{TYPEBAFF}
%and \ref{TYPECDAF} 
holds for $\widetilde{B}_3$.
\end{lemma}
{\em Proof.} Using $\imath_B$ (the transposition of
the two mirrors), one can assume that the first angle
of the configuration representing a minimal NGT,\,
$\{\tbe,\tga=\tal+\tbe,\tal\}$\,, 
is $\tal=\ep_2-\ep_3$. Then the
last one, $\tbe$, can be 

(1)\ $m\de +\ep_1+\ep_3$, or\ \  (2)\ $m\de +\ep_1-\ep_2$ 
\ for \ $m\ge 0$, and, additionally,

(3)\ $m\de -\ep_1-\ep_2$, or\ \ (4)\ $m\de -\ep_1+\ep_3$ 
\ when $m>0$.

Let us demonstrate that 
the last three choices are impossible. 
We will use Lemma \ref{BTINEQ}.

First of all, the following holds: \\
a) $L_2$ reflects in the bottom mirror
after the intersection with $L_3$,\\
b) the first reflection of $L_1$ may occur only in the 
bottom mirror, \\
c) $b_1\le b_2\ $ for the numbers of the bottom reflections
of $L_1$ and $L_2$,\\
d) the first reflection (if any)
of line $L_3$ can be only in the top mirror.

Furthermore, a simple check gives that the angles between 
$L_1$ and $L_2$ will be always in the form 
$m\de +\ep_2\pm \ep_1$; this excludes
(2) and (3). 

A more algebraic verification is as follows. 
If the angle from (2) for $m>0$ appears in the configuration,
then so does $\ep_1-\ep_2$.
The latter represents a simple root and can be made the 
first one, which contradicts the minimality of NGT. Similarly, 
for the angle from (3), $\de-\ep_1-\ep_2$ is an angle too; 
it results in a contradiction too.

As for (4), line $L_3$ intersects $L_1$ when it goes down
(after the corresponding top reflection) or up (after
the corresponding bottom reflection). In either case,
the sign of $\ep_1$ in the intersection angle 
is always plus, so (4) is impossible.
\medskip

\begin{figure}[htbp]
\begin{center}
\vskip 1.5in
\hskip -0.5in
\includegraphics[scale=0.45]{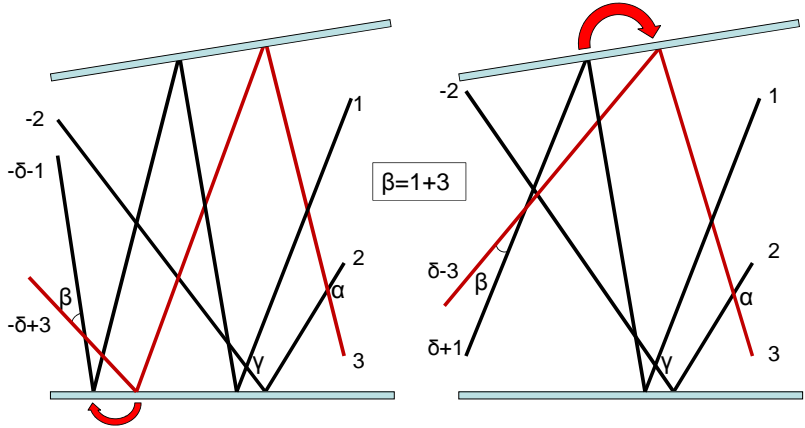}
\vskip -2.in
\caption{Line 3 must be horizontal}
\label{bngt2}
\end{center}
\end{figure}

Thus, (1) is the only option for $\tbe$. Let us now check
that line $L_3$ is actually horizontal (i.e., does not
reflect). We claim that if it reflects in the mirrors
then its last reflection can be made the last event of the
configuration, which contradicts the minimality of
the NGT under consideration.  Figure \ref{bngt2}
demonstrates this claim; the thick arc there shows the
reflection points that can be transposed in this
configuration.

%\vskip 1.5in
\medskip

This conclude the verification of the lemma.\sq 
\medskip

Let us apply Lemma  \ref{REDTO3} to the
three lines forming a minimal NGT
$\{\tbe,\tga,\tal\}$ for $\widetilde{B}_4$.
Then the forth line can be one of the following: \\
1) horizontal near the bottom or near the
top (the last or the first); \\
2) the one between the two non-horizontal 
parallel lines from the triple, \\
3) below these two non-horizontal lines and  
of type $b-t=1$ with a $b$\~number smaller
than that for these two.

Generally, for a minimal NGT configuration 
of type $\widetilde{B}_n$, any new line can be either 
added to an existing 
{\em bunch of lines} or can ``begin" a
new bunch subject to the inequalities 
from Section \ref{sec:B-NGT}.
It includes the horizontal bunches near the bottom or near the 
top. We use Lemma \ref{BTINEQ}. The theorem is proven.

\setcounter{equation}{0}
\section{Types
\texorpdfstring{\mathversion{bold}$C-D$}{\em C-D}} 
\begin{figure}[htbp]
\begin{center}
\vskip 1.5in
\hskip -0.5in
\includegraphics[scale=0.45]{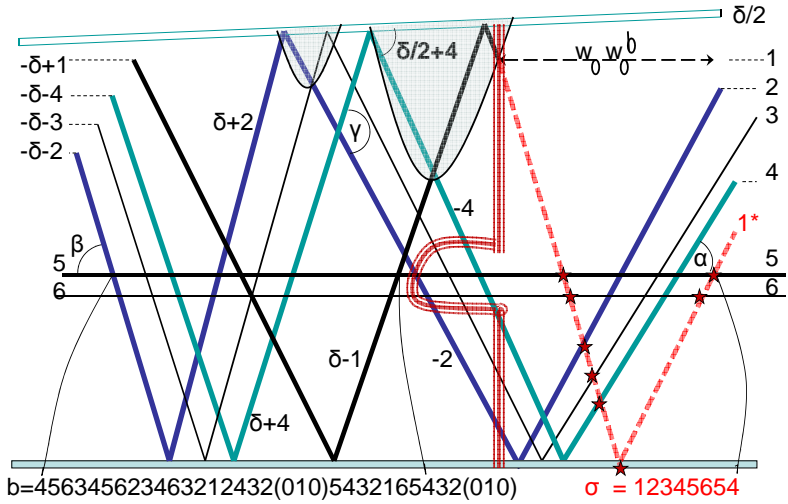}
\vskip -2.in
\caption{NGT of type C-D: $b$ and $\sigma$}
\label{cngt}
\end{center}
\end{figure}

Figure \ref{cngt} gives an example of a
minimal affine NGT of type $C$ constructed
using a {\em parity correction} from a
minimal NGT of type $B$.
It also illustrates the construction 
from Theorem \ref{ALMOSTPOSITIVE}, including
a geometric interpretation of $b$ and  $\si$. 

%\vskip 1.5in
\medskip

Let us use this graph
to demonstrate the changes in the $\widetilde{C},
\widetilde{D}$\~cases 
versus the planar interpretation for $\widetilde{B}$. 
We will use prime ($\tW'$, $\hW'$ and so on)
for the objects considered for $\widetilde{C}, \widetilde{D}$.
\smallskip

\subsection{Main modifications}
The general way of constructing the reduced decompositions 
in terms of the intersection and reflection points 
remains essentially the same. The elements $\tw\in \tW$ 
of type $\widetilde{B}$ always belong to $\hW'$ of 
type $\widetilde{C}$ or $\widetilde{D}$
(generally, not to $\tW'$); the element $s_0$ for $B$ is
naturally interpreted as the generator of $\Pi'$ for
$C$ or $D$. The following is necessary and sufficient
for the inclusion $\tw\in\tW'$. 

One needs to check that the total number of top reflections
is even (for $\widetilde{C}_n$ and  $\widetilde{D}_n$) 
and the total number of bottom reflections is even
in the case of $\widetilde{D}_n (n\ge 4)$.  
Then we can transform
such reduced decomposition to make it from $\tW'$,
i.e., in terms  of new simple events
$s_0'=s_0 s_1 s_0$ and $s_n'=s_n s_{n-1} s_n$; 
the latter is needed for $D$.

In Figure \ref{cngt}, if the %red 
dashed line to $1^*$
is added to line $1$, then the corresponding word becomes 
of type $\widetilde{D}_6$ (from the corresponding
$\tW$). If this dashed line
is disregarded here, then the corresponding configuration is
of type $\widetilde{C}$ but not of type $\widetilde{D}$.
If we disregard the dashed line completely and, moreover,
remove the top-right reflection of line $1$, then the
resulting word is neither of type $\widetilde{C}$ nor
of type $\widetilde{D}$. Let us discuss this
example and related features of our construction
in more detail.
\smallskip

First, let us begin with the element corresponding
to the configuration where we disregard
the portion of line $1$ {\em before} (to the right of)
the top-right reflection. Graphically, the dashed line
is disregarded. We will denote it by $\tw$.
Adding this top reflection to $\tw$ gives an example
of the {\em top-right parity correction} of $\tw$. 
Algebraically, $\tw\mapsto \tw'\equal\tw s_0$. It makes
the number of the top reflections even,
so $\tw'$ can be expressed in
terms of $s_0'$ instead of $s_0$ and becomes 
a word of type $\widetilde{C}$.

Second, let $\tw''$ be the graph for $\tw'$ extended by the %red 
dashed line ending at $1^{\star}$; then  
$\tw''$ is of type $D$.

\smallskip
Actually, the simplest way of transforming
the $\tw'$ to a word of type $D$ is via
the {\em bottom parity
correction} (right or left), i.e., using line $6$. 
Algebraically, it is the transformation
$\tw'\mapsto \tw's_n=s_n\tw'$. 

Note that the 
{\em top-left parity correction} of $\tw$, that is $s_0\tw$, 
is different from the {\em
top-right parity correction} $\tw'=\tw s_0$.
Generally, the right and left parity corrections coincide
only if they are performed on the same horizontal line.
Line $6$ (used for the 
bottom-right correction) is {\em horizontal};
line $1$ is not.  
\smallskip

Concerning the interpretation of the angles as roots
and related matters, there are the following modifications
versus the $\widetilde{B}$\~case. 
\medskip

($ i$) The angles for the bottom reflections 
must be multiplied by $2$ for $\widetilde{C}$.
The angle of 
$s'_0=s_0s_1s_0$ or $s_n'=s_{n}s_{n-1}s_{n}$,
presented in terms of $s_0,s_n$ for $\widetilde{B}$, is 
the middle one (from the three angles involved in this event).
\smallskip

($ ii$) The angles $j\de+p\pm q$,
including $j\de +2p$ in the $\widetilde{C}$\~case, are transformed
to the affine roots $[\,\ep_p\pm \ep_q\,,\,j\,]$; so 
this interpretation is different from the $\widetilde{B}$\~case, 
where the ``unit" was $\de/2$.
The graphic description of the action of $\tW$ on the roots 
remains unchanged; we read the angles after the event,
replacing their original numbers by the local ones.
 
In the figure under consideration, 
the angles of the two $\widetilde{D}$\~type top events (marked)
are correspondingly $\de-1+4=[\ep_1-\ep_4,1]$ and 
$\de+3+2=[\ep_2+\ep_3,1]$.
\smallskip

($ iii$) The interpretation of the sequence of
the {\em absolute angles} (with the $x$\~axis)
at the end of the configuration 
as a representation of $\tw$ in the form $bw$  
remains unchanged versus the $\widetilde{B}$\~case.

Recall that we consider the sequence of absolute
angles (counted from top to bottom) as a vector 
\begin{align}\label{debww}
&\de(b_1,b_2,\ldots,b_n)+w(1,2,\ldots,n)
\hbox{\ for proper }  b_i\in \Z \and 
w\in W.
\end{align}
Then $\tw=bw$, where we identify 
$b=(b_1,\ldots,b_2)$ with
$\sum_{i=1}^n b_i\ep_i\in Q$. See (\ref{debw}).
We continue using 
the notation from \cite{Bo}.  

For instance, $\al_0=[-\vth,1]$,
where $\vth=\ep_1+\ep_2$ for both, $\widetilde{C}$
and $\widetilde{D}$. The angle of $s_0'=s_0s_1s_0$ is 
$\de-1-2=[-\ep_1-\ep_2,1]$. The vector of the absolute angles
after this event is $\de(1,1)+(-2,-1)$. Thus
$b=\ep_1+\ep_2=\vth, w=s_{\vth}$ and 
$s_0'=\vth s_{\vth}$.
\medskip

Note that the lattice $Q$ becomes smaller versus that for 
$\widetilde{B}$
(it is the same one for $\widetilde{C}$ and $\widetilde{D}$). 
Namely,
it contains  $b=\sum_{i=1}^n b_i \ep_i$ with
integral $b_i$ only for even $\sum_{i=1}^n b_i$.

For instance, $\al_0=[-\vth,1]$,
where $\vth=\ep_1+\ep_2$ for both, $\widetilde{C}$
and $\widetilde{D}$. The angle of $s_0'=s_0s_1s_0\,$ is 
$\de-1-2=[-\ep_1-\ep_2,1]$. The vector of the absolute angles
after this event is $\de(1,1)+(-2,-1)$. Thus
$b=\ep_1+\ep_2=\vth$ and  $w=s_{\vth}$, which matches 
the relation $s_0'=\vth s_{\vth}$.

%\vskip 1.5in
\begin{figure}[htbp]
\begin{center}
\vskip 1.5in
\hskip -0.5in
\includegraphics[scale=0.45]{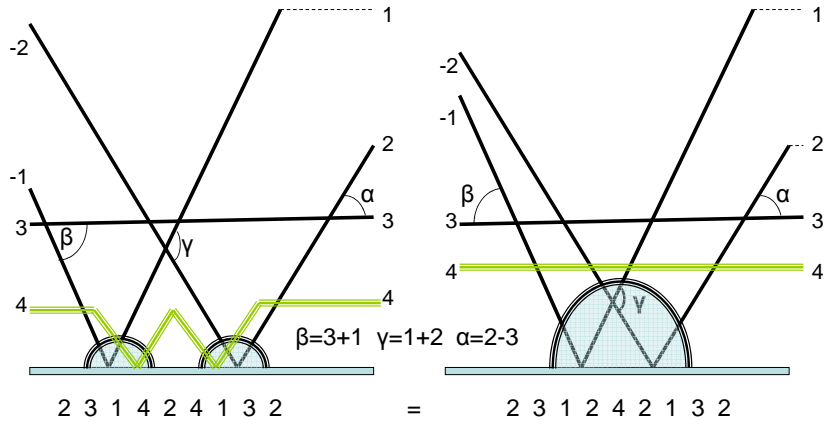}
\vskip -2.in
\caption{Type $D$, breaking the line}
\label{dngt}
\end{center}
\end{figure}
\medskip

\smallskip
Concerning $s_n'$ or $s_0'$, there is a special
procedure for dealing with the graphs when the other
lines are allowed to intersect (the area of)
the corresponding triple event.
It is called {\em breaking the line} and is directly
related to the parity corrections. It is necessary
for collecting the triples corresponding to $s_n'$ or $s_0'$
in a given $B$\~word (assuming that the latter satisfies 
the corresponding parity condition).

Figure \ref{dngt} reproduces the graph from
\cite{CS}, which demonstrates 
the procedure {\em breaking the line} and presents the simplest
nonaffine minimal NGT of type $D_4$. 
Here $4$ in the reduced decompositions
stays for the simple (nonaffine) reflection $s_4'$ for $D_4$,
corresponding to the triple product $s_4s_3s_4$
from the viewpoint of nonaffine $B$ or $C$. 
This graph proves the Coxeter 
relation $s'_{n}s_{n-2}s_n'=s_{n-2}s_n' s_{n-2}$
($424=242$).
\smallskip

\subsection{The classification theorem}
Theorem \ref{ALMOSTPOSITIVE}, generally, provides adequate
algebraic tools for the classical minimal NGT. 
As a matter of fact, 
its reduced version, Corollary \ref{SUBSYSNGT},
is sufficient to obtain all classical minimal
NGT if one uses the automorphisms of the affine
Dynkin diagram and the passage from $Q$ to $P$.
This corollary can be used to 
justify {\em algebraically} that our geometric
constructions (in the classical cases) really result in 
minimal NGT. However, concerning the completeness, we 
rely on the planar interpretation.
The algebraic approach based on Theorem \ref{ALMOSTPOSITIVE}
is expected to be useful in the theory of minimal
NGT for the exceptional root systems. 
\smallskip

The discussion above leads to the following theorem.
Slightly abusing the terminology, by a {\em mirror of type
$D$}, we mean the top mirror for  $\widetilde{C}$ or an 
either one for  $\widetilde{D}$ (i.e., when the 
corresponding end of the affine Dynkin diagram is of type $D$). 

\begin{theorem}
(i) The group $\hW'$ of type $\widetilde{C}$ can be naturally
identified with the group $\tW$ of type $\widetilde{B}$,
where  $s_0\in \tW$ is interpreted as the generator of 
$\Pi'\subset \hW'=\tW$. Accordingly, the configurations
for the elements from $\tW$  become configurations for
$\hW'$ and an arbitrary element of $\hW'$ can be obtained 
from a configuration of type $\widetilde{B}$.  

(ii) The resulting elements belong to $\tW'$ of type
 $\widetilde{C}$ if and only if the number of the top reflections 
of the corresponding configuration is even. Furthermore,
they belong to $\tW'$ of type $\widetilde{D}$ if the number
of bottom reflections is even too. Given $\tw'\in \tW'$,
it can be obtained from a $\widetilde{B}$\~type configuration 
where all $s_0$ and $s_n$ are included in the events 
$s_0'=s_0s_1s_0$ and $s_n'=s_ns_{n-1}s_n$. 

(iii) All minimal NGT from $\tW'$ of type  
$\widetilde{C}$ or $\widetilde{D}$ come from minimal
NGT of type $\widetilde{B}$ satisfying (ii).
However, the latter may become non-minimal in $\tW'$;
it occurs only if the horizontal line involved in the NGT is 
a unique one in its (horizontal) bunch and if the corresponding 
mirror is of type $D$.
\end{theorem}

{\em Proof}. The statements from ($i,ii$) have been 
already discussed. Concerning $\widetilde{C}$, we move all
$s_0$ to the beginning (or the end)  of a given reduced 
decomposition of
$\tw\in \tW$, replacing $s_0 s_1 s_0$ by $s_0'$ when 
necessary. It will give a word from $\tW'$ possibly
multiplied by $s_0$ on the right (or on the left).

As for $\widetilde{D}$, we can use that the group
$\tW$ for $\widetilde{B}$ can be naturally identified
with the extension of $\tW'$ of type $\widetilde{B}$
by $\Z_2\times\Z_2$ generated by the elements
$s_0,s_n$ (pairwise commutative) treated as outer automorphisms
of the corresponding affine Dynkin diagram. The element
$s_0$ is from $\hW'$, the element $s_n$ is not; both are
of zero length by definition. One can move all such
elements in a given reduced $\widetilde{B}$\~decomposition to its
beginning or to its end. The elements $s_0'$ and $s_n'$
may be produced during this process. The top (or bottom) 
parity corrections are needed if the elements 
$s_0$ (or $s_n$) do not cancel each other.

Algebraically, ($ii$) means that for any given $\tw'\in \tW'$,
its reduced decompositions with the minimal possible numbers
of $s_0'$ and $s_n'$ remain reduced in $\tW$,
where $s_0'$ and $s_n'$ are expressed in terms of $s_0$ and
$s_n$. The geometric approach guarantees that at least one
such reduced decomposition exists. The construction
$\tw\mapsto\tw'$ (subject to the parity conditions) consists
of moving all $s_0,s_n$ to the beginning (or to the end) of a
given reduced decomposition of $\tw$. 
\smallskip

Concerning ($iii$), 
let us use the right graph in Figure \ref{dngt} 
to demonstrate what
happens in the beginning of the configuration
if there is only one line parallel to the corresponding
mirror. If (the lowest) line $4$ is removed, then the
corresponding reduced decomposition reads
$s_2s_1s_3's_1s_2$. Using the Coxeter $D$\~relation,
it can be transformed to 
$s_2s_3's_1s_3's_2=s_3's_2s_1s_2s_3'$. Thus, the beginning
of the resulting reduce decomposition is movable (not unique)
and such $B$\~minimal NGT will not remain $D$\~minimal.

This example is actually a general one; 
it is sufficient to manage arbitrary $\widetilde{C}_n$ and
$\widetilde{D}_n$ if the horizontal line involved in
the NGT is near the mirror of type $D$ and is the only one in
its horizontal bunch. Then the right end of the resulting 
$\la$\~sequence becomes movable upon the switch to 
$\widetilde{C}_n$ or
$\widetilde{D}_n$ in this case (and only in such case). 
\sq
\smallskip

For instance, let us consider the configuration from 
Figure \ref{cngt} including the dashed line to $1^{\star}$ 
and excluding line $6$. It 
represents a minimal NGT for $\widetilde{B}_5$,
which {\em will not} remain minimal upon its
recalculation
to $\widetilde{D}_5$ (which is possible because the 
numbers of top and bottom reflections are both even).
It is analogous 
to the constraint ``at least two horizontal bottom lines" 
from \cite{CS} in the nonaffine $D$\~case. 
\smallskip

The following theorem is an explicit form of claim ($iii$),
reformulated in terms of the parity corrections.  

\begin{theorem}\label{TYPECDAF}
(i) An arbitrary minimal NGT from $\tW$
of type $\widetilde{C}_n (n\ge 3)$ can be obtained
as $\tw'$, $\tw''$ or $\tw^*$ as follows.

Let $\tw$ be a $B$\~positive minimal NGT of type
$\widetilde{B}_n (n\ge 3)$. If no top
parity correction is needed for $\tw$ 
(i.e., the total number of top reflections in $\tw'$ 
is even), then $\tw'\equal \tw$, 
(considered as $\widetilde{C}$\~words) and, moreover,
$\tw''=\imath_C(\tw)=s_0\tw s_0$ are minimal NGT of type 
$\widetilde{C}_n$. The element $\tw^*$ is
minimal NGT if it has even number of
the top reflections and the initial
$\tw$ has at least two bottom 
horizontal lines; $\tw^*$ is involutive and 
coincides with $\imath_C(\tw^*)$.
 
Otherwise, if the parity corrections are needed,
$\tw'\equal\tw s_0$ (the results of the top-right parity 
correction) is a minimal NGT of $\widetilde{C}_n$\~type. 
For the top-left correction, $\tw''\equal s_0\tw=
\imath_C(\tw')=(\tw')^{-1}$, so such element can be
represented as $\tu'$ for certain $B$\~positive $\tu$.
The element $\tw^*$ defined now as 
$\imath_B(\tw) s_0=s_0 \imath_B(\tw)$, 
is a minimal NGT subject to the same 
``2 line constraint" as above. 

(ii) Minimal NGT of type $D$ are given in
terms of the $C$\~words $\tw',\tw'',\tw^*$ from part 
(i) as follows. If no bottom parity correction is needed, 
i.e., the total number of the bottom reflections
in the initial $\tw$ (the cases of $\tw'$ and $\tw''$)
or $\imath_B(\tw)$ (the case of $\tw^*$) is even, 
then each of these elements is a minimal NGT of type 
$\widetilde{D}_n$. We require here the existence of at least
two horizontal lines near the bottom for $\tw$
in the cases of $\tw'$ and $\tw''$. 

Otherwise, if the total number
of the bottom reflections for $\tw'$, $\tw''$ or
for $\tw^*$ is odd, then the elements 
$s_n\tw'=\tw' s_n$ and $s_n\tw''=\tw'' s_n$, as well
as the elements 
$\tw^*s_n$ and  $s_n\tw^*=(\tw^* s_n)^{-1}$
are minimal NGT of $D$\~type. We
impose the same constraint as above in the cases 
of $\tw',\tw''$, namely, the configuration 
for $\tw$ is supposed to contain at least two 
horizontal lines near the bottom.
\sq  
\end{theorem}  

We note that the operation $\tw\mapsto s_0\tw s_0$
for $B$\~positive $\tw$ is trivial if
$\tw$ contains at least one {\em top horizontal}
line. It is obvious geometrically that
the configuration for $s_0\tw s_0$ can be transformed
to ensure the cancelation of such two $s_0$ if 
they are ``performed" on the same top horizontal line.
Similarly, the transformation $\tw\mapsto
s_n \tw s_n$ is not needed because it is
always trivial; the bottom line of 
$B$\~positive $\tw$ is always {\em horizontal}.
\medskip

\subsection{Algebraic aspects}
Let us now discuss Figure \ref{cngt} from the
viewpoint of ``algebraic" Theorem \ref{ALMOSTPOSITIVE} and 
Corollary \ref{SUBSYSNGT}.
The $C$\~elements
$\tw'$ and $\tw''$ (the latter is $\tw'$
extended by the %red
dashed line to $1^{\star}$) have
coinciding weights, so they must correspond to different $\si$.
This weight is $b=-(\ep_1+\ep_2+\ep_3+\ep_4)$
in the notation from \cite{Bo};
the portion of the graph after (to the left of) 
the solid thick %brown 
vertical ``line" gives its reduced 
decomposition. The ``enclave" there excludes the intersections 
with lines $5$ and $6$ near this thick line, which simply 
means that we are supposed to make these two lines 
``very" close to the bottom to make the decomposition
of $b$ right.

We note that 
$l(\tw')=l(b)+l(\upsilon_b^\si)$
for $\tw'=b \upsilon_b^\si$ from the figure, which is
generally not the case. It occurs if and only if
$(b,\al)\le 0$ for all $\al\in R_+$.

For the configuration for $tw''$
(that includes the dashed line to $1^{\star}$),
$\si$ is identical, i.e.,
$(\upsilon_b^\si)^{-1}=(w_0 w_0^b)^{-1}$ here.
The element $w_0\dot{w}_0^b$ (recall that $w_0=-1$)
is given by the portion of this configuration 
before (to the right of) the thick %brown 
``line".  
It is shown as $w_0w_0^b$ in this picture. 
The {\em dot}\~marks can be omitted here and below
because $p>0$ and $b$ is {\em not} dominant, namely, 
$p=-(b,\al_4)=-(b,\ep_4-\ep_5)=1.$

The example of $\tw'=\varpi_b^\si$ (with $1$ instead of 
$1^{\star}$), 
requires using a non-trivial $\si$ (see the graph), 
resulting in ``deleting" the dashed line.
Its reduced decomposition is $12345654$ (where we 
put $i$ instead of $s_1$).

From the viewpoint of Corollary \ref{SUBSYSNGT},
we consider the extension from $\tR'=\widetilde{B}_5$
to $\tR=\widetilde{C}_6$, where  
$b'=-(\ep_2+\ep_3+\ep_4)$ is extended to
$b=b'-\ep_1$. Then $\tw'=b w_0'w_0^{b'}$,
where $w_0'w_0^{b'}$ is calculated in
$\tW'$ for $\widetilde{B}_5$. It gives
exactly $b(w_0\si w_0^b)^{-1}$.
\smallskip

\rmk
Note, that the corresponding $\varpi_b^\si$ is not 
a minimal NGT when considered as a word of type $B$.
It is obvious from the picture.
Algebraically, $\si(\vth)=\vth=\ep_1$ in 
$B_6$ and the corresponding condition from 
(\ref{vthcond}),  $(b,\vth)\le 1$, fails.
Recall, that our form is normalized by the
condition $(\vth,\vth)=2$, so $(b,\vth)=2$.

When doing the same check for $\varpi$ treated
as a $C$\~word, we use that $\vth$ becomes $\ep_1+\ep_2$ 
and our form
coincides with the standard one (such that
$\{\ep_j\}$ are orthonormal). Then
$(b,\vth)=0\le 0$, which guarantees that
the outcome is a minimal NGT.  \sq
\smallskip

%\vskip 1.5in
\begin{figure}[htbp]
\begin{center}
\vskip 1.5in
\hskip -0.5in
\includegraphics[scale=0.45]{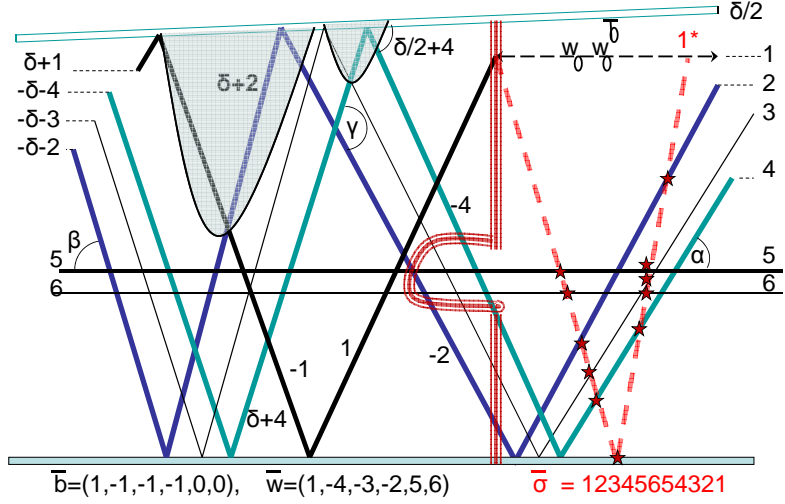}
\vskip -2.in
\caption{Type $C$: top-left correction}
\label{cngt1}
\end{center}
\end{figure}
\medskip

{\bf The inversion.}
We will begin with an analysis of the
inversion of the minimal NGT constructed
above. 
Continuing the discussion of 
Theorem \ref{ALMOSTPOSITIVE}, let us consider 
$\tw^\#\equal(\varpi_b^\si)^{-1}$ equals 
$\varpi_{\bar{b}}^{\bar{\si}}$, where
$\bar{b}=\si(\varsigma(b))=\ep_1-\ep_2-\ep_3-\ep_3$.
Recall that $\tw'=\varpi_b^\si$; see Figure \ref{cngt}.
Graphically, $\tw^\#$ describes
the {\em top-left correction} of the element 
$\tw$ above. 

Figure \ref{cngt1} shows the configuration 
associated with  $\tw^\#=\bar{b}\bar{w}$;
it corresponds to the case when the
dashed line is ignored.
Geometrically, it is quite obvious that
it represents a minimal NGT.
Let us interpret this fact algebraically
via Theorem \ref{ALMOSTPOSITIVE}. 

The weights $b$ and
$\bar{b}$ belong to the same set $Q\cap P_+^{(4)}$
(have non-negative inner products with the same $\al_4$),
i.e., $\bar{k}=k=4$.
The condition ($ a'$) from (\ref{upsilonsharp}) (which
guarantees that $\al_{k}$ is a non-movable first root
in $\la(\varpi_b^\si$) holds. Moreover, ($ a'''$) 
from (\ref{lengdbc}) holds. Namely,
$\si(R_{+}^{\,b})\supset R_+^{\bar{b}}$.
Indeed, 
$$\si(\al_1,\al_2,\al_3,\al_5,\al_6)=
(\al_2,\al_3,\al_4,\al_5,\al_6)\supset 
(\al_2,\al_3,\al_5,\al_6). 
$$
For the sake of completeness, let us mention
that $\si(\al_4)=-\ep_1-\ep_5$.

The element $\bar{\si}$ equals $12345654321$; so it 
is non-trivial and does
not coincides with $\si^{-1}$. In Figure \ref{cngt1},
$\bar{\si}$ is represented by the portion of the configuration
before (to the right of) the vertical line 
where only the events involving the dashed line to $1^*$ are
taken.

The element $\bar{b}$ is shown by the
portion after (to the left of) the vertical line. Notice
that the events for $s_0'$ are combined in a way
different from that in Figure \ref{cngt}. The element
$w_0w_0^{\bar{b}}$ corresponds to the portion of
the configuration taken before the vertical line
and where the dashed line to $1*$ is included; 
it is shown in the picture. 

From the viewpoint of Corollary \ref{SUBSYSNGT},
we consider the extension from $\tR'=\widetilde{B}_5$
to $\tR=\widetilde{C}_6$, where  
$b'=-(\ep_2+\ep_3+\ep_4)$ is extended to
$\bar{b}=b'+\ep_1$. Then $\tw^{\#}=\bar{b} w_0'w_0^{b'}$,
where $w_0'w_0^{b'}$ is calculated in
$\tW'$ for $\widetilde{B}_5$. 

The condition ($ b'''$) from (\ref{lengdbc})
is not satisfied for $\bar{b}$.
However, ($b''$) from 
(\ref{lengdbcnewb}) holds:
$$
\bar{\si}(R_{+}^{\,\bar{b}})\,\cup\, \la(\bar{\si}^{-1}) 
\supset R_{+}^{b},
$$ 
where the root $\al_1\in  R_{+}^{b}$, missing in
$\bar{\si}(R_{+}^{\,\bar{b}})$, comes 
from $\la(\bar{\si}^{-1})$.
\smallskip

We note that the total word in
Figure \ref{cngt1},
{\em including the dashed line}, 
is {\em non-reduced}. 
It represents the element $\tw^{\flat}\equal\varpi_{\bar{b}}$ 
with the trivial $\si$. 
This element is a minimal NGT of type
$\widetilde{C}$ and also of type 
$\widetilde{D}$. 
The corresponding reduced word
can be obtained  by transforming line $1$ (including the dashed 
portion) into a horizontal one; its two bottom reflections 
annihilate each other. The resulting word for $\tw^\flat$
is a result of the  {\em top parity correction}
of a $B$\~positive minimal NGT in a {\em horizontal}
line, namely, line $1$ upon making it horizontal.
\medskip

%\vskip 1.5in
\begin{figure}[htbp]
\begin{center}
\vskip 1.5in
\hskip -0.5in
\includegraphics[scale=0.45]{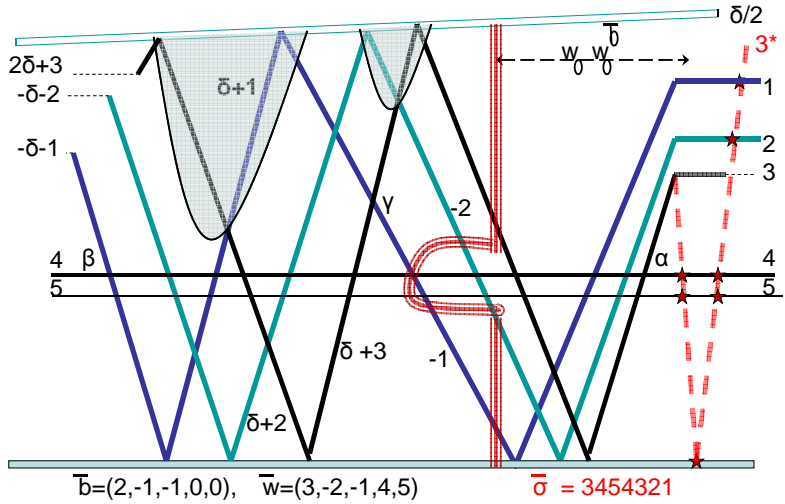}
\vskip -2.in
\caption{Type $C$: non-trivial $\si$}
\label{cngt2}
\end{center}
\end{figure}
\medskip
 
{\bf More on the parity corrections.}
The most involved ones are those performed on a line
that is already involved in the {\em affine} reflections
in the initial minimal NGT. Corollary 
\ref{SUBSYSNGT} covers such parity corrections
only if we apply it for a complete $\hW$ (with
$P$ instead of $Q$) with further reduction to $\tW$
(the division by $\pi_r$). See (\ref{pitrans}).

A typical example of this situation is as follows.
Let $\tw'=b'w'$ for $b'=(-1,-1,-1,0,0)$ and 
$w'=(-3,-2,-1,0,0)$.
In terms of the data $(a,b,c)$, this $B$\~positive
element can be described as follows. First,
$u=0$ (there is no top horizontal bunch of lines).
Second, there is only one bunch from ($b$) with $t_1=1$, 
which contains $3$ lines.
Third, the bottom horizontal bunch contains $2$
lines ($v=2$). The corresponding element $\tw$
is of type $\widetilde{B}_5$ but not of type
$\widetilde{C}_5$, since it involves $3$ top reflections.

Let us perform the top-right correction. The resulting
element $\tw$ is of type $\widetilde{C}$
(but not of type $\widetilde{D}$). One has:
$$
\tw=bw \for 
b=(-1,-1,-2,0,0) \and w=(-3,-2,1,4,5).
$$ 
The element $w$ here is not of type $w_0w_0^c$ 
(including the modifications from  
Corollary \ref{SUBSYSNGT}). 
See (\ref{d-example}) for the algebraic
meaning of the parity correction we used; 
the corresponding weight there is denoted by $b''$.
We omit the graph of $\tw$. Instead, we will
focus on $\tw^*=\tw^{-1}$. 

Figure \ref{cngt2} shows the
configuration for $\tw^*$; the dashed line must be
dropped:
$$
\tw^*=\bar{b}\bar{w} \for 
\bar{b}=(2,-1,-1,0,0) \and \bar{w}=(3,-2,1,4,5).
$$ 
This element is a result of the top-left correction
of the initial element $\tw'$.
The expression for $\bar{\si}$ is provided in the
figure. See (\ref{d-example}); the corresponding
weight there is $b'$.

The algebraic mechanism ensuring the absence
of $\al_0$ in $\la(\tw)$ and $\la(\tw^*)$ 
is of some interest.
The inner product $(\bar{b},\vth')=1$ is critical
for the second condition from formula 
(\ref{vthcond}); here $\vth'=\ep_1+\ep_2$.
So we need to check that 
$\bar{\si}w_0^{\bar{b}}(\vth')>0$. One has:
$$
w_0^{\bar{b}}(\vth')=\ep_1+\ep_3,\ 
 \bar{\si}(\ep_1+\ep_3)=\ep_2-\ep_3>0.
$$
 
If the dashed line is added to $\tw^*$, then the
corresponding element $\tw^\#$
is a result of the two
simultaneous top parity corrections of the
$B$\~word $\tw''$ constructed for $b''=(0,-1,-1,0,0)$
and $w''=(-1,-3,-2,4,5)$. See Figure \ref{cngt2}:
$\tw^\#=s_0 \tw'' s_0$.
The corresponding $\si$\~factor is trivial for
$\tw^\#$, so it is of the simplest type
from the viewpoint of 
Theorem \ref{ALMOSTPOSITIVE}. 

%An interpretation
%of $\tw, \tw^*$ and similar ones entirely in terms of 
%$\widetilde{C}$ (without using the parity corrections)
%is unclear at the moment. Such clarification 
%is expected to lead to counterparts of the
%elements of this type for exceptional root systems. Any
%other types of classical minimal NGT can be generalized.   
\medskip

\subsection*{Acknowledgements}
The paper is partially supported by NSF grant DMS--0800642.
The first author is thankful to Maxim Nazarov for
interesting discussions on twisted Yangians 
(and his recent results with Sergey Khoroshkin)
and for the invitation to the York University.
The major part of the paper was written 
during the stay of I.Ch. at the University Paris 7;
he is grateful to Eric Vasserot and the 
Foundation for Mathematical Sciences of Paris
for the invitation.

%\vskip -1cm
%\medskip
\bibliographystyle{unsrt}

\end{document}